\documentclass[journal,transmag]{IEEEtran}

\usepackage{cite}
\usepackage{amsmath}
\usepackage{booktabs}
\usepackage{tabularx}
\usepackage{subcaption}
\usepackage{import}
\usepackage{bm}
\usepackage{amsfonts}
\usepackage{amssymb}
\usepackage{tikz}
\usepackage{soul}

\usetikzlibrary{babel,shapes.geometric,arrows,calc,positioning}
\newcommand\myeq{\mathrel{\overset{\makebox[0pt]{\mbox{\normalfont\tiny 2-D}}}{=}}}
\newcommand\marc[1]{{\color{black}{#1}}} 

\begin{document}

\title{Hyperbolic Basis Functions for Time-Transient Analysis of Eddy Currents in Conductive and Magnetic Thin Sheets}

\author{\IEEEauthorblockN{Bruno de Sousa Alves\IEEEauthorrefmark{1},
Ruth V. Sabariego\IEEEauthorrefmark{2},
Marc Laforest\IEEEauthorrefmark{1}, and
Fr\'ed\'eric Sirois\IEEEauthorrefmark{1}}
\IEEEauthorblockA{\IEEEauthorrefmark{1} Polytechnique Montr\'eal, Montr\'eal, QC, Canada}
\IEEEauthorblockA{\IEEEauthorrefmark{2} Department of Electrical Engineering (ESAT), EnergyVille, KU Leuven, 3001 Leuven, Belgium}
}

\IEEEtitleabstractindextext{
\begin{abstract}

This paper presents a new time-domain finite-element approach for modelling thin sheets with hyperbolic basis functions derived from the well-known steady-state solution of the linear flux diffusion equation. The combination of solutions at different operating frequencies permits the representation of the time-evolution of field quantities in the magnetic field formulation. This approach is here applied to solve a planar shielding problem in harmonic and time-dependent simulations for materials with either linear or nonlinear characteristics. Local and global quantities show good agreement with \textcolor{black}{the} reference solutions obtained by the standard finite element method on a complete and representative discretization of the region exposed to a time-varying magnetic field.
\end{abstract}

}

\maketitle
\IEEEdisplaynontitleabstractindextext
\IEEEpeerreviewmaketitle

\section{Introduction}
\IEEEPARstart{T}{hin} sheets of high permeability and/or conductivity are often employed to mitigate stray fields produced by electric and electronic devices such as rotating machines, large power transformers, induction heating equipment, welding and forming machines~\cite{bottauscio2006transient,Igarashi1998a,Igarashi1998b,bottauscio2004numerical,Rasilo2020}. In the surroundings of these devices, the field intensity needs to be at acceptable value to comply with \textcolor{black}{the} exposure limits for humans and for electromagnetic compatibility and interference reasons~\cite{bottauscio2006transient}. However, the shielding efficiency of thin sheets can be directly affected by its material characteristics, shape and position~\cite{bottauscio2004numerical}. Consequently, the availability of models able to predict the electromagnetic (EM) behavior in such structures at an affordable computational cost is key to optimizing these devices.

In terms of shape, the high aspect ratio of thin-sheet shields presents a challenge to numerical simulation.  Indeed, the direct application of a numerical method, such as the Finite Element Method (FEM), can be computationally expensive or even prohibitive due to the associated meshing difficulties~\cite{geuzaine2000dual}. On the one hand, a coarse mesh inside these sheets is unable to capture the EM phenomena and may lead to elements with high aspect ratio, which affect the FEM solution accuracy and convergence~\cite{marchandise2014optimal}. On the other hand, a high-density mesh can increase inordinately the number of unknowns in the problem and therefore the computational cost. The EM problem is even harder to solve if nonlinear materials characteristics are \textcolor{black}{considered} in time-transient analysis.   

An efficient way to overcome these difficulties is to use the classical Thin-Shell (TS) \textcolor{black}{\mbox{model~\cite{krahenbuhl1993thin,mayergoyz1995calculation,Guerin1994,Biro1997}}.} In this model, a reduced-dimension geometry replaces the actual volume of the thin regions, and suitable \textit{impedance boundary conditions}~{(IBCs)} account for the EM behavior \marc{within the original volume}.  These IBCs are defined from the analytical solution of the EM problem throughout the volume of the thin sheet. Thus, since the \marc{smallest dimension} of the layer is neglected in a geometric sense, \marc{errors are avoided that would have been} caused by the \marc{original anisotropic} meshing of the thin structure \marc{with poor aspect ratios}~\cite{Igarashi1998a}.

The TS model provides a good compromise between accuracy and computational cost~\cite{geuzaine2000dual}, but its application is still mostly restricted to linear and harmonic regime analysis since the analytical solution is known \textit{a~priori} and the IBCs can be \textcolor{black}{easily} established.  Currently available time-domain and nonlinear TS-FEM approaches are often derived from classical IBCs, whether using orthogonal polynomial basis functions to express the magnetic flux density through the shell \mbox{thickness~\cite{gyselinck2008,sabariego2008h,sabariego2009,gyselinck2004time,gyselinck2006,gyselinck2015finite},} Fast Fourier Transform~(FFT) to update the residual from the harmonic solution~\cite{bottauscio2004numerical,bottauscio2006transient}, or simply a linear field variation through its thickness (strip approximations)~\cite{brandt1996superconductors,Zhang2017,Liang2017,Berrospe2019}. However, when considering time-transient analysis of nonlinear thin sheets, \textcolor{black}{more representative models are required.}

In~\cite{Sabariego2010}, low-order \textit{surface impedance boundary conditions} (SIBCs) are defined using basis functions derived from the steady-state analytical solution of semi-infinite slab problems. These SIBCs are applied in time-transient FE simulations to remove large conducting regions from the computational domain. Although the problems involving thin sheets are different, their nature is the same.  \textcolor{black}{In SIBCs, field quantities penetrate the surface from one boundary of a bulk domain, whereas in the TS model, the penetration occurs simultaneously from the two extended faces of a thin sheet.} Moreover, the IBCs in the TS approach proposed in~\cite{mayergoyz1995calculation} are derived in a similar way than the SIBCs in~\cite{mayergoyz1994finite}. For these reasons, a time-domain extension of the classical TS model with basis functions derived from the steady-state solution of a slab of finite thickness, \textcolor{black}{equivalent to the model proposed in~\cite{Sabariego2010} for SIBCs,} is a natural approach to pursue.  

This paper presents a novel time-domain extension of the classical TS model to solve 2-D shielding problems. The physics inside the thin region is \textcolor{black}{captured} by hyperbolic basis functions derived from the steady-state analytical solution of the linear flux diffusion equation. We demonstrate that the use of \textcolor{black}{two hyperbolic basis functions} leads to IBCs equivalent to the classical TS model in harmonic regime. In the time-transient analysis, the use of $n$ \textcolor{black}{pairs of hyperbolic} basis functions, \marc{each representing different frequencies, coupled to FEM models outside the TS,} allows computation of the time evolution of the physical quantities \marc{throughout the domain without resolving the TS.} The method  is here developed for a magnetic field ($\bm{h}$-)formulation and extended to nonlinear cases. Results show good agreement with the 2-D FE reference solutions, with a greatly reduced number of degrees of freedom (DoFs) and therefore, at a lower computational cost. 

\section{1-D Flux Diffusion Problem in a Slab of Finite Thickness}

\par
\begin{figure}[t]
	\centering
	\begin{subfigure}{.2\textwidth}
		\centering
		\tikzset{every picture/.style={line width=0.4pt}} 

\begin{tikzpicture}[x=0.6pt,y=0.6pt,yscale=-0.8,xscale=0.8]

\draw  [draw opacity=0][fill={rgb, 255:red, 200; green, 200; blue, 200 }  ,fill opacity=1 ] (72,120.25) -- (236,120.25) -- (236,170.5) -- (72,170.5) -- cycle ;
\draw    (72,144.5) -- (114.75,144.62) ;
\draw [shift={(116.75,144.63)}, rotate = 180.16] [color={rgb, 255:red, 0; green, 0; blue, 0 }  ][line width=0.75]    (10.93,-3.29) .. controls (6.95,-1.4) and (3.31,-0.3) .. (0,0) .. controls (3.31,0.3) and (6.95,1.4) .. (10.93,3.29)   ;
\draw    (72,144.5) -- (72.24,92.12) ;
\draw [shift={(72.25,90.13)}, rotate = 450.26] [color={rgb, 255:red, 0; green, 0; blue, 0 }  ][line width=0.75]    (10.93,-3.29) .. controls (6.95,-1.4) and (3.31,-0.3) .. (0,0) .. controls (3.31,0.3) and (6.95,1.4) .. (10.93,3.29)   ;
\draw    (122,113.75) -- (153,113.98) ;
\draw [shift={(156,114)}, rotate = 180.42] [fill={rgb, 255:red, 0; green, 0; blue, 0 }  ][line width=0.08]  [draw opacity=0] (8.93,-4.29) -- (0,0) -- (8.93,4.29) -- cycle    ;
\draw    (122,177.75) -- (153,177.98) ;
\draw [shift={(156,178)}, rotate = 180.42] [fill={rgb, 255:red, 0; green, 0; blue, 0 }  ][line width=0.08]  [draw opacity=0] (8.93,-4.29) -- (0,0) -- (8.93,4.29) -- cycle    ;
\draw    (192,119.75) -- (192,100.88) ;
\draw [shift={(192,97.88)}, rotate = 450] [fill={rgb, 255:red, 0; green, 0; blue, 0 }  ][line width=0.08]  [draw opacity=0] (8.93,-4.29) -- (0,0) -- (8.93,4.29) -- cycle    ;
\draw    (192,189.75) -- (192,170.88) ;
\draw [shift={(192,192.75)}, rotate = 270] [fill={rgb, 255:red, 0; green, 0; blue, 0 }  ][line width=0.08]  [draw opacity=0] (8.93,-4.29) -- (0,0) -- (8.93,4.29) -- cycle    ;
\draw    (72,120.25) -- (236,120.25) ;
\draw    (72,170.5) -- (236,170.5) ;
\draw  [fill={rgb, 255:red, 255; green, 255; blue, 255 }  ,fill opacity=1 ] (68.05,144.07) .. controls (68.05,141.75) and (69.93,139.87) .. (72.25,139.87) .. controls (74.57,139.87) and (76.46,141.75) .. (76.46,144.07) .. controls (76.46,146.39) and (74.57,148.27) .. (72.25,148.27) .. controls (69.93,148.27) and (68.05,146.39) .. (68.05,144.07) -- cycle ;
\draw    (269,120.25) -- (269,170.5) ;
\draw    (265.5,120.08) -- (272.64,119.94) ;
\draw    (265.5,170.08) -- (272.64,169.94) ;
\draw  [fill={rgb, 255:red, 0; green, 0; blue, 0 }  ,fill opacity=1 ] (71.15,144.07) .. controls (71.15,143.46) and (71.64,142.96) .. (72.25,142.96) .. controls (72.86,142.96) and (73.36,143.46) .. (73.36,144.07) .. controls (73.36,144.68) and (72.86,145.18) .. (72.25,145.18) .. controls (71.64,145.18) and (71.15,144.68) .. (71.15,144.07) -- cycle ;

\draw (104.75,151) node [anchor=north west][inner sep=0.75pt]    {$x$};
\draw (75.75,89.75) node [anchor=north west][inner sep=0.75pt]    {$y$};
\draw (55.61,130.96) node [anchor=north west][inner sep=0.75pt]    {$z$};
\draw (123.25,86.75) node [anchor=north west][inner sep=0.75pt]    {$h^{+}_{x}$};
\draw (123.75,178) node [anchor=north west][inner sep=0.75pt]    {$h^{-}_{x}$};
\draw (238,106.5) node [anchor=north west][inner sep=0.75pt]    {$\Gamma ^{+}_{s}$};
\draw (236.5,155.5) node [anchor=north west][inner sep=0.75pt]    {$\Gamma ^{-}_{s}$};
\draw (195.5,87) node [anchor=north west][inner sep=0.75pt]    {$\bm{n}^{+}_{s}$};
\draw (195.25,172.25) node [anchor=north west][inner sep=0.75pt]    {$\bm{n}^{-}_{s}$};
\draw (193,135) node [anchor=north west][inner sep=0.75pt]    {$\Omega _{s}$};
\draw (270.14,131.71) node [anchor=north west][inner sep=0.75pt]    {$d$};

\end{tikzpicture}
		\vspace{-3mm}
    	\caption{Magnetic field} 
    	\label{Thin_regionH}
	\end{subfigure}
	\hspace{0.6cm}
	\begin{subfigure}{.2\textwidth}
		\centering
		\tikzset{every picture/.style={line width=0.4pt}} 

\begin{tikzpicture}[x=0.6pt,y=0.6pt,yscale=-0.8,xscale=0.8]

\draw  [draw opacity=0][fill={rgb, 255:red, 200; green, 200; blue, 200 }  ,fill opacity=1 ] (315.8,120.25) -- (479.8,120.25) -- (479.8,170.5) -- (315.8,170.5) -- cycle ;
\draw    (315.8,144.5) -- (358.55,144.62) ;
\draw [shift={(360.55,144.63)}, rotate = 180.16] [color={rgb, 255:red, 0; green, 0; blue, 0 }  ][line width=0.75]    (10.93,-3.29) .. controls (6.95,-1.4) and (3.31,-0.3) .. (0,0) .. controls (3.31,0.3) and (6.95,1.4) .. (10.93,3.29)   ;
\draw    (315.8,144.5) -- (316.04,92.12) ;
\draw [shift={(316.05,90.13)}, rotate = 450.26] [color={rgb, 255:red, 0; green, 0; blue, 0 }  ][line width=0.75]    (10.93,-3.29) .. controls (6.95,-1.4) and (3.31,-0.3) .. (0,0) .. controls (3.31,0.3) and (6.95,1.4) .. (10.93,3.29)   ;
\draw    (365.8,113.75) -- (396.8,113.98) ;
\draw [shift={(399.8,114)}, rotate = 180.42] [fill={rgb, 255:red, 0; green, 0; blue, 0 }  ][line width=0.08]  [draw opacity=0] (8.93,-4.29) -- (0,0) -- (8.93,4.29) -- cycle    ;
\draw    (365.8,177.75) -- (396.8,177.98) ;
\draw [shift={(399.8,178)}, rotate = 180.42] [fill={rgb, 255:red, 0; green, 0; blue, 0 }  ][line width=0.08]  [draw opacity=0] (8.93,-4.29) -- (0,0) -- (8.93,4.29) -- cycle    ;
\draw    (435.8,119.75) -- (435.8,100.88) ;
\draw [shift={(435.8,97.88)}, rotate = 450] [fill={rgb, 255:red, 0; green, 0; blue, 0 }  ][line width=0.08]  [draw opacity=0] (8.93,-4.29) -- (0,0) -- (8.93,4.29) -- cycle    ;
\draw    (435.8,189.75) -- (435.8,170.88) ;
\draw [shift={(435.8,192.75)}, rotate = 270] [fill={rgb, 255:red, 0; green, 0; blue, 0 }  ][line width=0.08]  [draw opacity=0] (8.93,-4.29) -- (0,0) -- (8.93,4.29) -- cycle    ;
\draw    (315.8,120.25) -- (479.8,120.25) ;
\draw    (315.8,170.5) -- (479.8,170.5) ;
\draw  [fill={rgb, 255:red, 0; green, 0; blue, 0 }  ,fill opacity=1 ] (315.59,144.53) .. controls (315.59,144.44) and (315.67,144.36) .. (315.77,144.36) .. controls (315.86,144.36) and (315.94,144.44) .. (315.94,144.53) .. controls (315.94,144.63) and (315.86,144.71) .. (315.77,144.71) .. controls (315.67,144.71) and (315.59,144.63) .. (315.59,144.53) -- cycle ;
\draw    (512.8,120.25) -- (512.8,170.5) ;
\draw    (509.3,120.08) -- (516.44,119.94) ;
\draw    (509.3,170.08) -- (516.44,169.94) ;
\draw  [fill={rgb, 255:red, 255; green, 255; blue, 255 }  ,fill opacity=1 ] (311.56,144.36) .. controls (311.56,142.04) and (313.45,140.15) .. (315.77,140.15) .. controls (318.09,140.15) and (319.97,142.04) .. (319.97,144.36) .. controls (319.97,146.68) and (318.09,148.56) .. (315.77,148.56) .. controls (313.45,148.56) and (311.56,146.68) .. (311.56,144.36) -- cycle ;
\draw   (312.63,141.53) .. controls (314.33,139.82) and (317.04,139.77) .. (318.68,141.4) .. controls (320.31,143.04) and (320.26,145.74) .. (318.55,147.45) .. controls (316.85,149.15) and (314.15,149.2) .. (312.51,147.57) .. controls (310.87,145.93) and (310.93,143.23) .. (312.63,141.53) -- cycle ; \draw   (312.63,141.53) -- (318.55,147.45) ; \draw   (318.68,141.4) -- (312.51,147.57) ;

\draw (348.55,151) node [anchor=north west][inner sep=0.75pt]    {$z$};
\draw (319.55,89.75) node [anchor=north west][inner sep=0.75pt]    {$y$};
\draw (297.07,131.63) node [anchor=north west][inner sep=0.75pt]    {$x$};
\draw (367.05,86.42) node [anchor=north west][inner sep=0.75pt]    {$e^{+}_{z}$};
\draw (367.55,178) node [anchor=north west][inner sep=0.75pt]    {$e^{-}_{z}$};
\draw (481.8,106.5) node [anchor=north west][inner sep=0.75pt]    {$\Gamma ^{+}_{s}$};
\draw (480.3,155.5) node [anchor=north west][inner sep=0.75pt]    {$\Gamma ^{-}_{s}$};
\draw (439.3,87) node [anchor=north west][inner sep=0.75pt]    {$\bm{n}^{+}_{s}$};
\draw (439.05,172.25) node [anchor=north west][inner sep=0.75pt]    {$\bm{n}^{-}_{s}$};
\draw (436.8,135) node [anchor=north west][inner sep=0.75pt]    {$\Omega _{s}$};
\draw (513.94,131.71) node [anchor=north west][inner sep=0.75pt]    {$d$};

\end{tikzpicture}
		\vspace{-3mm}
		\caption{Electric field}
		\label{Thin_regionE}
	\end{subfigure}
\vspace{-2mm}
	\caption{\small Thin region $\Omega_s$ \textcolor{black}{of thickness $d$} and its local coordinate system. \marc{$\Gamma_s^{\pm}$
	denote the top and bottom boundaries of $\Omega_s$ and $\boldsymbol{n}_s^{\pm}$ their respective outward normals.}}
	\label{Thin_region}
	\vspace{-4mm}
\end{figure}
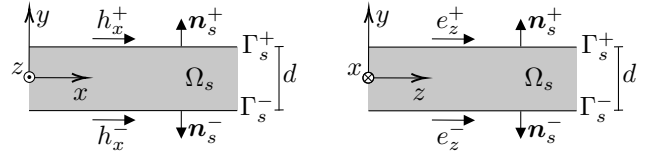

 By assuming a thin region with a sufficiently high aspect ratio, the EM problem in a sheet can be formulated as a \mbox{1-D} flux diffusion problem in a slab of finite thickness. In Fig.~\ref{Thin_region}, we consider a thin sheet of thickness $d$ whose normal is parallel to the \mbox{$y$-axis}. The tangential component of the magnetic field ($h_x$) is in the \mbox{$x$-direction} \mbox{(Fig.~\ref{Thin_regionH})} and the tangential component of the electric field (${e}_z$) is in the \mbox{$z$-direction} \mbox{(Fig.~\ref{Thin_regionE})}. The slab problem can be then  formulated in terms of ${h}_x$ or ${e}_z$, i.e.
\begin{align}
\partial_y(\rho \partial_y {h}_x(y,t) ) + \partial_t \mu {h}_x(y,t) =0, \label{1D_h} \\
\partial_y (\nu  \partial_y {e}_z(y,t) )+ \partial_t \sigma {e}_z(y,t) =0, \label{1D_e} 
\end{align} 
where $\rho$ is the electric resistivity ($\sigma=1/\rho$) and $\mu$ is the magnetic permeability ($\nu=1/\mu$). These expressions are derived from Faraday's and Ampere's laws, respectively. 

We assume that $\rho$ and $\mu$ are constants and that we have harmonic boundary conditions~(BCs), i.e.
\begin{align}
{h}_{x}^\pm(y=\pm d/2,t) = \hat{h}_x^\pm \cos(\omega t +\phi_h^\pm), \label{BC1D_h} \\
{e}_{z}^\pm(y=\pm d/2,t) = \hat{e}_z^\pm \cos(\omega t +\phi_e^\pm), \label{BC1D_e} 
\end{align} 
where $\hat{h}_x^\pm$ and $\hat{e}_z^\pm$ are respectively the magnetic and electric field magnitudes, and  $\phi_h^\pm$ and $\phi_e^\pm$ their related phase shifts. Using a phasor representation (symbols with a bar), i.e. \textcolor{black}{\mbox{$\bar{h}_{x}^\pm =\hat{h}_x^\pm \exp({\jmath \phi_h^\pm})$} and \mbox{$\bar{e}_{z}^\pm =\hat{e}_z^\pm \exp({\jmath \phi_e^\pm})$,}} we have
\begin{align}
{h}_{x}^\pm(y=\pm d/2,t)= \Re\{ \bar{h}_{x}^\pm \exp({\jmath \omega t})\},   \label{BC1D_h_complex} \\
{e}_{z}^\pm(y=\pm d/2,t)= \Re\{ \bar{e}_{z}^\pm \exp({\jmath \omega t})\},  \label{BC1D_e_complex} 
\end{align}
\textcolor{black}{where $\Re\{.\}$ is the real part of the argument.}

\marc{The solutions to~\eqref{1D_h}-\eqref{1D_e} with the boundary \mbox{conditions~\eqref{BC1D_h_complex}-\eqref{BC1D_e_complex}} are given complex notation} (symbols with right arrow on top) and can be written as:
\begin{align}
\bar{h}_x(y) = \bar{h}_{x}^+\vec{\psi}^+ (y) + \bar{h}_{x}^-\vec{\psi}^- (y), \label{Sol1D_h} \\
\bar{e}_z(y) = \bar{e}_{z}^+\vec{\psi}^+ (y) + \bar{e}_{z}^-\vec{\psi}^- (y), \label{Sol1D_e}
\end{align}
where 
\begin{equation} \label{psi_basis}
\vec{\psi}^\pm (y)=\frac{\sinh \left(\frac{\vec{a} d}{2} \pm \vec{a} y\right)}{\sinh\left(\vec{a} d\right)},
\end{equation}
and $\vec{a}=\frac{1+\jmath}{\delta}$, \textcolor{black}{$\jmath=\sqrt{-1}$,}  $\delta = \sqrt{2/(\mu\sigma\omega)}$, $\omega=2\pi f$, and $f$ is the operating frequency. Note that the functions $\vec{\psi}^\pm(y)$ appear in both~\eqref{Sol1D_h} and~\eqref{Sol1D_e}. These functions are used later in this paper to define the basis functions required for the time-transient analysis of thin sheets.

\subsection{Hyperbolic Basis Functions and Classical TS Model}

We propose the use of basis functions obtained from steady-state solutions of the 1-D flux diffusion problem \textcolor{black}{governed} by~\eqref{1D_h}-\eqref{BC1D_e}. Taking as example the problem in terms of the magnetic field~\eqref{1D_h}, together with the BCs~\eqref{BC1D_h}, the steady-state solution for $h_x$ can be written as
\vspace{-2mm}
\begin{equation} \label{sincos}
\begin{split}
{h}_{x}(y,t) = \,&\hat{h}_x^+ \cos(\omega t +\phi_h^+) \theta_c^+(y) \\
&+ \hat{h}_x^+ \sin(\omega t + \phi_h^+) \theta_s^+(y)\\
&+  \hat{h}_x^- \cos(\omega t + \phi_h^-) \theta_c^-(y)\\
&+  \hat{h}_x^- \sin(\omega t + \phi_h^+) \theta_s^-(y),
\end{split}
\vspace{-1mm}
\end{equation}
where \textcolor{black}{the analytical expressions for $\theta_c^\pm$ and $\theta_s^\pm$ given in Table~4.2-II of~\cite{knoepfel2008magnetic}.} \textcolor{black}{Here,} $\theta_c$ and $\theta_s$ are obtained \textcolor{black}{directly} \textcolor{black}{from the real~($\Re$) and imaginary~($\Im$) parts of~\eqref{psi_basis},} i.e.
\vspace{-1mm}
\begin{align}
\theta_c^\pm(y)=\Re\{\vec{\psi}^\pm(y)\}, \\
\theta_s^\pm(y)=\Im\{\vec{\psi}^\pm(y)\},
\end{align} \vspace{-8mm}

\noindent which means that
\vspace{-2mm}
\begin{equation}  \label{psitheta}
\vec{\psi}^\pm(y)  = \theta_c^\pm(y)+\jmath \theta_s^\pm(y).
\vspace{-2mm}
\end{equation}

\begin{figure} [t]
	\hspace{-0.3cm}
	\includegraphics[width=0.5\textwidth]{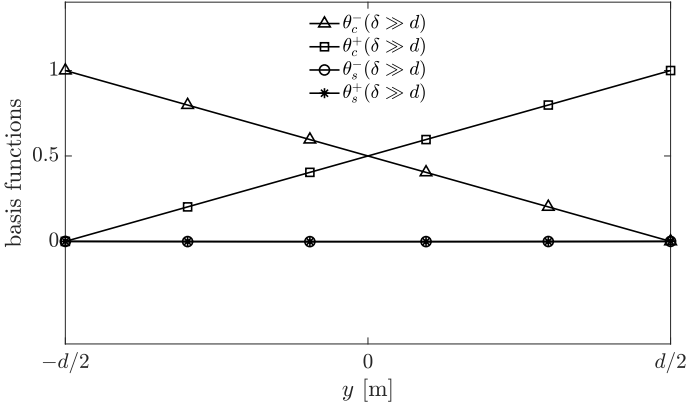}
	\caption{\small Hyperbolic basis functions when $\delta \gg d$.}
	\vspace{-0.5cm}
	\label{Basis_deltaggd}
\end{figure}

Note that the solution~\eqref{sincos} can be interpreted as a least squares approximation of $h_x^\pm$ in $\Omega_s$ using the hyperbolic functions $\vec{\psi}^\pm(y)$. In addition, if we consider the harmonic solution~\eqref{Sol1D_h} with $\delta \gg d$, $\vec{\psi}^\pm(y)$ \textcolor{black}{in~\eqref{psi_basis}} reduces to
\begin{equation} \label{psi}
\left. \vec{\psi}^\pm(y) \right|_{(\delta \gg d)}= \frac{d/2 \pm y}{d},
\end{equation} 
which is equivalent to the Lagrange polynomials of first order defined across the thickness $d$ of the sheet. Indeed, with $\delta \gg d$, the field quantities have a linear variation through the sheet thickness, and the functions $\vec{\psi}^\pm(y)$ can account for this behavior. To illustrate this, the functions $\theta_c^\pm$ and $\theta_s^\pm$ for $\delta \gg d$ are plotted in Fig.~\ref{Basis_deltaggd}.

The equivalent solution in terms of the electric field can be obtained by \marc{replacing} $\hat{h}_x^\pm$ and $\phi_h^\pm$ by $\hat{e}_z^\pm$ and $\phi_e^\pm$ in~\eqref{sincos}. However, in this paper, we are solely interested in magnetic field quantities and the $\bm{h}$-formulation. The solution depending on the electric field would be useful\marc{, say, to implement} the proposed approach in the magnetic vector potential \mbox{($\bm{a}$-)formulation.} 

Besides, it can be demonstrated that the application of the \textcolor{black}{$\vec{\psi}^\pm(y)$} functions as basis functions in the variational form of a \mbox{1-D} finite element problem in the harmonic regime leads to \marc{the same} IBCs used in the classical TS model~\cite{mayergoyz1995calculation}, which are, \textcolor{black}{using a vector representation (bold symbols),}
\vspace{-1mm} 
\begin{align}
{\bm{n}}_s \times ({\bm{h}}_{x}^+ - {\bm{h}}_{x}^-) & = {\vec{\eta}}_e \left( {\bm{n}}_s \times \left( {\bm{e}}_{z}^+ + {\bm{e}}_{z}^-\right)\right) \times {\bm{n}}_s, \label{IBC01} \\
{\bm{n}}_s \times ({\bm{e}}_{z}^+ - {\bm{e}}_{z}^-) & = {\vec{\eta}}_h \left( {\bm{n}}_s \times \left( {\bm{h}}_{x}^+ + {\bm{h}}_{x}^-\right)\right) \times {\bm{n}}_s, \label{IBC02}
\end{align}
with ${\vec{\eta}}_h = -\frac{\jmath\omega \mu}{{\vec{a}}} \tanh \left(\frac{{\vec{a}} d}{2}\right)$ and ${\vec{\eta}}_e = \frac{\sigma}{{\vec{a}}} \tanh \left(\frac{{\vec{a}} d}{2}\right)$. 

Expression~\eqref{IBC01} connects the discontinuity of the tangential components of the magnetic field to the mean value of the tangential electric field. This discontinuity is related to the total net current flowing in the sheet~\cite{krahenbuhl1993thin}. Moreover, equation~\eqref{IBC02} connects the discontinuity of the tangential electric field to the mean value of the tangential magnetic field, which is related to the amount of perpendicular flux absorbed in the plane of the sheet. When $\delta \gg d$, the coefficients  $\frac{1}{\vec{a}}\tanh \left(\frac{{\vec{a}} d}{2}\right)$ in ${\vec{\eta}}_h$ and ${\vec{\eta}}_e$ can be approximated by $d/2$~\cite{geuzaine2000dual}. 

\textcolor{black}{The classical TS model  in the form of IBCs has been used extensively to tackle problems involving thin regions in harmonic regime simulations\cite{mayergoyz1995calculation,krahenbuhl1993thin,geuzaine2000dual,Guerin1994,bottauscio2004numerical,Rasilo2020}. These IBCs were originally defined from the analytical solution for the field distribution, and the integration of the analytical expressions of the electric and magnetic current densities over the thickness of the thin sheet, which gives equivalent surface currents representatives~\cite{mayergoyz1995calculation}. However, to the best of our knowledge, no definition of the TS model in the form of the hyperbolic basis functions has been proposed in the literature before. \textcolor{black}{These functions appear naturally in the solution~\eqref{sincos}.} 
\textcolor{black}{Therefore, the proposed approach} can be easily extended to time-transient analysis, as described next.}

\subsection{Hyperbolic Basis Functions in Time-Transient Analysis} \label{Dedicated}

\par In time-transient analysis, we define $n$ pairs of~$\vec{\psi}_k^\pm(y)$, where $k$ is the harmonic rank relative to a fundamental frequency $f_1$ chosen in accordance with the problem to model, and $1\leq k \leq n$. The number of basis functions $n$ is defined according to the frequency content of $h_x$ and the desired accuracy. Therefore, $\vec{\psi}_k^\pm(y)$ is still defined by~\eqref{psi_basis}, with $\omega=2\pi f_k$ (which affects the values of $\delta$ and $\vec{a}$). Then, according to~\eqref{psitheta}, each $\vec{\psi}^\pm_k(y)$ generates the even~$\theta_{ck}^\pm(y)$ and odd~$\theta_{sk}^\pm(y)$ functions, which we write as
\vspace{-1mm}
\begin{align}
\theta_{c1}^\pm(y) &= \Re\{\vec{\psi}_1^\pm(y)\}, \label{Basis1} \\
\theta_{ck}^\pm(y) &=  \Re\{\vec{\psi}_k^\pm(y)\} - \theta_{c1}^\pm(y), \hspace{1.3cm} \ 2\leq k \leq n, \label{Basis2} \\
\theta_{sk}^\pm(y) &=  \Im\{\vec{\psi}_k^\pm(y)\}, \hspace{2.8cm} 1\leq k \leq n. \label{Basis3}
\vspace{-1mm}
\end{align}

\begin{figure} [t]
	\hspace{-0.3cm}
	\includegraphics[width=0.5\textwidth]{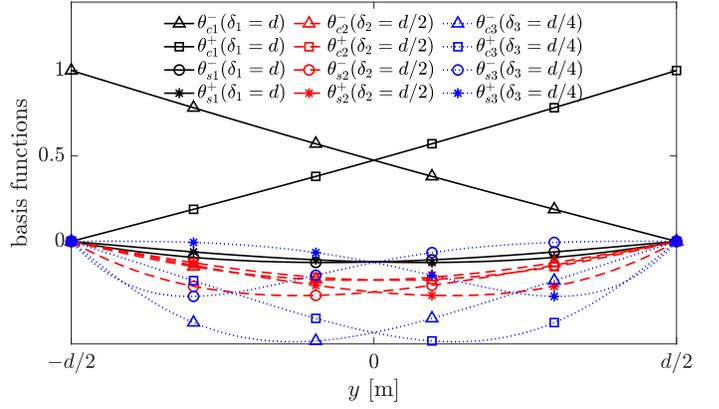}
	\caption{\small Hyperbolic basis functions for time-transient analysis: example with $\delta_1 = d$ (with $f_1$ in accordance), $\delta_2 = d/2$ ($f_2=4f_1$) and $\delta_3 = d/4$ ($f_3=16f_1$)  for $k=1,2$ and 3, respectively.}
	\label{Basis_deltasd3}
	\vspace{-3.7mm}
\end{figure}

\marc{The first two cosines satisfy \mbox{$\theta_{c1}^\pm(y=\pm d/2)=1$} while the remaining functions $\theta_{ck}^\pm(y)$ and $\theta_{sk}^\pm(y)$ in~(\ref{Basis2}-\ref{Basis3}) vanish at the boundaries of the thin region ($\Gamma^\pm$)}. This allows us to connect the \mbox{1-D} equations to the exterior FE global system of equations; see Section~\ref{FEMimple}. Examples of the proposed basis functions with $n=3$ are presented in Fig.~\ref{Basis_deltasd3} for $\delta \leq d$. For cases with $\delta \gg d$, the first two cosines functions are enough to represent the profile of $h_x$ in $\Omega_s$, since it has a linear variation throughout the thickness of the sheet (Fig.~\ref{Basis_deltaggd}). \marc{Furthermore, the sine} basis vanish everywhere.

\par The expansion of $h_x(y,t)$ in terms of~\eqref{Basis1}-\eqref{Basis3} can be written in matrix form as
\vspace{-1mm}
\begin{equation} \label{H_inside}
h_x(y,t) = \left[h(t)\right]^T \left[\theta(y)\right],
\end{equation}
\vspace{-1mm}
with the $4n\times 1$ matrices $\left[h(t)\right]$ and $\left[\theta(y)\right]$ given by
\vspace{-0.5mm}
\begin{equation}
\begin{split}
\left[h(t)\right] = [h_{c1}^+(t) \dots h_{cn}^+(t), 
h_{s1}^+(t) \dots h_{sn}^+(t),\\
h_{c1}^-(t) \dots h_{cn}^-(t),
h_{s1}^-(t) \dots h_{sn}^-(t)],
\end{split}
\end{equation}
\vspace{-1mm}
\begin{equation} \label{thetacs}
\begin{split}
\left[\theta(y)\right] = [\theta_{c1}^+(y) \dots \theta_{cn}^+(y),
\theta_{s1}^+(y) \dots \theta_{sn}^+(y),\\
\theta_{c1}^-(y) \dots \theta_{cn}^-(y),
\theta_{s1}^-(y) \dots \theta_{sn}^-(y)  ],
\end{split}
\end{equation}\vspace{-1mm}\hspace{-2mm}where $h_{ck}^\pm(t)$ and $h_{sk}^\pm(t)$ are unknowns of the problem to model.

Additionally, the 1-D variational form of the partial differential equation~\eqref{1D_h}, \textcolor{black}{disregarding homogeneous BCs,} is
\vspace{-1mm}
\begin{equation} \label{h_weak1D}
\Big(\rho \; \partial_y h_x, \partial_y h_x' \Big)_{\Omega_s} + \partial_t \Big(\mu \; h_x, h_x' \Big)_{\Omega_s} = 0,
\end{equation}\vspace{-0.5mm}where $h_x'$ is the test function \marc{assumed to vanish at $\Gamma_s^{\pm}$}.
\par The FE discretization of~\eqref{h_weak1D} by means of $N=4n$ basis functions $\theta_p(y)$ and $\theta_q(y)$, with $p,q\in [1,N]$, for $h_x$ and $h_x'$ respectively, and assuming isotropic linear materials, leads to a system of equations, expressed in matrix form as
\begin{equation} \label{EqSys}
\rho\left[\mathcal{S}\right] \left[h(t)\right] + \mu\left[ \mathcal{M} \right]\partial_t[h(t)] = 0,
\end{equation}
where the elements of $[\mathcal{S}]$ and $[\mathcal{M}]$ are calculated as
\begin{equation} \label{Spq}
\mathcal{S}_{pq} = \int_{-d/2}^{d/2} \partial_y \theta_p(y)\partial_y\theta_q(y) dy,
\end{equation}
\begin{equation} \label{Mpq}
\mathcal{M}_{pq} = \int_{-d/2}^{d/2} \theta_p(y) \theta_q(y) dy,
\end{equation}
which can be evaluated numerically for each pair of basis functions $\theta_p$ and $\theta_q$. Then, considering the implicit Euler scheme for the time-discretization of~\eqref{EqSys}, coupled to~\eqref{BC1D_h}, one obtains a system of algebraic equations to be solved at each time-step of the simulation. 

The instantaneous loss $\mathcal{L}(t)$ in \textcolor{black}{Joule} is calculated as~\cite{Gyselinck2009}
\begin{equation}
\mathcal{L}(t) = \rho [h(t)]^T [\mathcal{S}] [h(t)],
\end{equation}
where $\mathcal{S}_{pq}$ is given \textcolor{black}{by~\eqref{Spq}.}

In the nonlinear case, the still isotropic resistivity~$\rho$ (or the magnetic permeability $\mu$) in the variational form~\eqref{h_weak1D} can depend on the magnetic field intensity~$h_x$ or its derivative~$\partial_yh_x$. The resulting nonlinear system of equations is solved by the Newton-Raphson (NR) iterative method, as presented in~\cite{Sabariego2010}, but with integral terms evaluated over the thickness of the thin region, i.e., $-d/2\leq y \leq d/2$. These integrals are solved numerically using the Legendre-Gauss quadrature at every iteration of the NR method. 

\section{FEM Implementation} \label{FEMimple}

\begin{figure}[t]
	\centering
	\begin{subfigure}{0.2\textwidth}
		\centering
		\tikzset{every picture/.style={line width=0.4pt}} 

\begin{tikzpicture}[x=0.6pt,y=0.6pt,yscale=-0.7,xscale=0.7]
\hspace{-5mm}

\draw  [fill={rgb, 255:red, 200; green, 200; blue, 200 }  ,fill opacity=1 ] (250.83,138.03) .. controls (260.06,121.75) and (277.86,119.04) .. (283.79,139.38) .. controls (289.73,159.73) and (310.17,163.8) .. (295.66,183.47) .. controls (281.16,203.13) and (258.74,206.52) .. (243.57,181.43) .. controls (228.41,156.34) and (241.59,154.3) .. (250.83,138.03) -- cycle ;
\draw    (204.74,76.8) -- (194.4,57.27) ;
\draw [shift={(193,54.62)}, rotate = 422.11] [fill={rgb, 255:red, 0; green, 0; blue, 0 }  ][line width=0.08]  [draw opacity=0] (8.93,-4.29) -- (0,0) -- (8.93,4.29) -- cycle    ;
\draw   (71,97.95) .. controls (88.55,74.21) and (134.27,83.08) .. (177,80.95) .. controls (210.67,80.12) and (202,72.95) .. (233,65.95) .. controls (251.02,63.73) and (262.88,70.47) .. (273,76.95) .. controls (312,101.95) and (332,158.95) .. (334,203.95) .. controls (336,248.95) and (290,263.95) .. (244,251.95) .. controls (198,239.95) and (183,241.95) .. (140,255.95) .. controls (97,269.95) and (74,269.95) .. (57,248.45) .. controls (40,226.95) and (37,143.95) .. (71,97.95) -- cycle ;
\draw  [fill={rgb, 255:red, 155; green, 155; blue, 155 }  ,fill opacity=1 ] (132.4,215.65) .. controls (124.97,214.81) and (118.71,211.87) .. (114.53,206.69) .. controls (102.26,191.49) and (112.62,162.79) .. (137.66,142.58) .. controls (162.7,122.38) and (192.94,118.32) .. (205.21,133.52) .. controls (208.96,138.17) and (210.6,144.09) .. (210.34,150.63) -- (200.77,154.33) .. controls (201.79,148.48) and (200.91,143.29) .. (197.83,139.48) .. controls (188.85,128.35) and (164.58,133.04) .. (143.61,149.96) .. controls (122.65,166.88) and (112.93,189.61) .. (121.91,200.74) .. controls (125.32,204.97) and (130.94,206.91) .. (137.76,206.77) -- cycle ;
\draw    (93.2,205.2) -- (109.45,196.17) ;
\draw [shift={(111.2,195.2)}, rotate = 510.95] [color={rgb, 255:red, 0; green, 0; blue, 0 }  ][line width=0.75]    (10.93,-3.29) .. controls (6.95,-1.4) and (3.31,-0.3) .. (0,0) .. controls (3.31,0.3) and (6.95,1.4) .. (10.93,3.29)   ;
\draw    (167.6,135.67) -- (158.74,116.08) ;
\draw [shift={(157.5,113.35)}, rotate = 425.65] [fill={rgb, 255:red, 0; green, 0; blue, 0 }  ][line width=0.08]  [draw opacity=0] (8.93,-4.29) -- (0,0) -- (8.93,4.29) -- cycle    ;
\draw    (175.11,124.63) -- (182.35,144.04) ;
\draw [shift={(183.4,146.85)}, rotate = 249.54000000000002] [fill={rgb, 255:red, 0; green, 0; blue, 0 }  ][line width=0.08]  [draw opacity=0] (8.93,-4.29) -- (0,0) -- (8.93,4.29) -- cycle    ;
\draw    (130.21,148.55) -- (119.95,136.53) ;
\draw [shift={(118,134.25)}, rotate = 409.51] [fill={rgb, 255:red, 0; green, 0; blue, 0 }  ][line width=0.08]  [draw opacity=0] (8.93,-4.29) -- (0,0) -- (8.93,4.29) -- cycle    ;
\draw    (137.47,156.32) -- (147.75,167.26) ;
\draw [shift={(149.8,169.45)}, rotate = 226.79] [fill={rgb, 255:red, 0; green, 0; blue, 0 }  ][line width=0.08]  [draw opacity=0] (8.93,-4.29) -- (0,0) -- (8.93,4.29) -- cycle    ;
\draw    (92.96,162.18) -- (108.02,172.55) ;
\draw [shift={(109.67,173.69)}, rotate = 214.56] [color={rgb, 255:red, 0; green, 0; blue, 0 }  ][line width=0.75]    (10.93,-3.29) .. controls (6.95,-1.4) and (3.31,-0.3) .. (0,0) .. controls (3.31,0.3) and (6.95,1.4) .. (10.93,3.29)   ;
\draw    (138.55,195.04) -- (125.23,184.68) ;
\draw [shift={(123.65,183.45)}, rotate = 397.86] [color={rgb, 255:red, 0; green, 0; blue, 0 }  ][line width=0.75]    (10.93,-3.29) .. controls (6.95,-1.4) and (3.31,-0.3) .. (0,0) .. controls (3.31,0.3) and (6.95,1.4) .. (10.93,3.29)   ;

\draw (147.82,260) node [anchor=north west][inner sep=0.75pt]    {$\Gamma =\Gamma _{h} \cup \Gamma _{e}$};
\draw (254.97,152.14) node [anchor=north west][inner sep=0.75pt]    {\small $\Omega _{c}$};
\draw (199.2,39.06) node [anchor=north west][inner sep=0.75pt]    {\small ${\bm{n}}_{\Gamma }$};
\draw (78.24,223) node [anchor=north west][inner sep=0.75pt]    {$\Omega _{c}^{C}$};
\draw (71,193.84) node [anchor=north west][inner sep=0.75pt]  [rotate=-0.29]  {\small $\Omega _{s}$};
\draw (146.78,87.5) node [anchor=north west][inner sep=0.75pt]  [rotate=-0.29]  {\small ${\bm{n}}_{\Gamma }^{-}$};
\draw (170,138.54) node [anchor=north west][inner sep=0.75pt]  [rotate=-0.29]  {\small ${\bm{n}}_{\Gamma }^{+}$};
\draw (150,158) node [anchor=north west][inner sep=0.75pt]  [rotate=-0.29]  {\small ${\bm{n}}_{s}^{-}$};
\draw (103.47,108.99) node [anchor=north west][inner sep=0.75pt]  [rotate=-0.29]  {\small ${\bm{n}}_{s}^{+}$};
\draw (139.32,183.18) node [anchor=north west][inner sep=0.75pt]  [rotate=-0.29]  {\small $\Gamma _{s}^{-}$};
\draw (68,146) node [anchor=north west][inner sep=0.75pt]  [rotate=-0.29]  {\small $\Gamma _{s}^{+}$};

\end{tikzpicture}
		\vspace{-2mm}
		\caption{}
		\vspace{-2mm}
		\label{fig:domainFull}
	\end{subfigure} \hspace{0.6cm}
	\begin{subfigure}{0.2\textwidth}
		\tikzset{every picture/.style={line width=0.4pt}} 

\begin{tikzpicture}[x=0.6pt,y=0.6pt,yscale=-0.7,xscale=0.7]
\hspace{-5mm}

\draw  [fill={rgb, 255:red, 200; green, 200; blue, 200 }  ,fill opacity=1 ] (565.83,134.95) .. controls (575.06,118.68) and (592.86,115.96) .. (598.79,136.31) .. controls (604.73,156.65) and (625.17,160.72) .. (610.66,180.39) .. controls (596.16,200.06) and (573.74,203.45) .. (558.57,178.36) .. controls (543.41,153.26) and (556.59,151.23) .. (565.83,134.95) -- cycle ;
\draw    (522.74,72.92) -- (512.4,53.39) ;
\draw [shift={(511,50.74)}, rotate = 422.11] [fill={rgb, 255:red, 0; green, 0; blue, 0 }  ][line width=0.08]  [draw opacity=0] (8.93,-4.29) -- (0,0) -- (8.93,4.29) -- cycle    ;
\draw   (389,94.08) .. controls (406.55,70.33) and (452.27,79.2) .. (495,77.08) .. controls (528.67,76.24) and (520,69.08) .. (551,62.08) .. controls (569.02,59.85) and (580.88,66.59) .. (591,73.08) .. controls (630,98.08) and (650,155.08) .. (652,200.08) .. controls (654,245.08) and (608,260.08) .. (562,248.08) .. controls (516,236.08) and (501,238.08) .. (458,252.08) .. controls (415,266.08) and (392,266.08) .. (375,244.58) .. controls (358,223.08) and (355,140.08) .. (389,94.08) -- cycle ;
\draw  [fill={rgb, 255:red, 155; green, 155; blue, 155 }  ,fill opacity=1 ] (450.85,203.57) .. controls (444.06,202.38) and (438.29,199.35) .. (434.28,194.38) .. controls (421.87,179.01) and (431.11,150.97) .. (454.91,131.77) .. controls (478.71,112.57) and (508.06,109.46) .. (520.46,124.84) .. controls (524.73,130.13) and (526.44,136.92) .. (525.87,144.33) -- (525.87,144.33) .. controls (526.44,136.92) and (524.73,130.13) .. (520.46,124.84) .. controls (508.06,109.46) and (478.71,112.57) .. (454.91,131.77) .. controls (431.11,150.97) and (421.87,179.01) .. (434.28,194.38) .. controls (438.29,199.35) and (444.06,202.38) .. (450.85,203.57) -- cycle ;
\draw    (480.57,116.82) -- (475.91,102.26) ;
\draw [shift={(475,99.4)}, rotate = 432.28] [fill={rgb, 255:red, 0; green, 0; blue, 0 }  ][line width=0.08]  [draw opacity=0] (8.93,-4.29) -- (0,0) -- (8.93,4.29) -- cycle    ;
\draw    (421.23,128) -- (437.49,138.64) ;
\draw [shift={(439.17,139.73)}, rotate = 213.19] [color={rgb, 255:red, 0; green, 0; blue, 0 }  ][line width=0.75]    (10.93,-3.29) .. controls (6.95,-1.4) and (3.31,-0.3) .. (0,0) .. controls (3.31,0.3) and (6.95,1.4) .. (10.93,3.29)   ;
\draw    (465.32,157.65) -- (449,146.18) ;
\draw [shift={(447.37,145.03)}, rotate = 395.11] [color={rgb, 255:red, 0; green, 0; blue, 0 }  ][line width=0.75]    (10.93,-3.29) .. controls (6.95,-1.4) and (3.31,-0.3) .. (0,0) .. controls (3.31,0.3) and (6.95,1.4) .. (10.93,3.29)   ;

\draw (569.97,149.06) node [anchor=north west][inner sep=0.75pt]    {$\Omega _{c}$};
\draw (465.82,256.12) node [anchor=north west][inner sep=0.75pt]    {$\Gamma =\Gamma _{h} \cup \Gamma _{e}$};
\draw (391.24,213.07) node [anchor=north west][inner sep=0.75pt]    {$\Omega _{c}^{C}$};
\draw (517.2,35.19) node [anchor=north west][inner sep=0.75pt]    {${\bm{n}}_{\Gamma }$};
\draw (475,84) node [anchor=north west][inner sep=0.75pt]  [rotate=-0.29]  {\small ${\bm{n}}_{s}$};
\draw (465.49,146.8) node [anchor=north west][inner sep=0.75pt]  [rotate=-0.29]  {\small $\Gamma _{s}^{-}$};
\draw (392.39,107.44) node [anchor=north west][inner sep=0.75pt]  [rotate=-0.29]  {\small $\Gamma _{s}^{+}$};
%

\draw  [fill={rgb, 255:red, 0; green, 0; blue, 0 }  ,fill opacity=1 ] (522.52,144.33) .. controls (522.52,142.98) and (523.57,141.89) .. (524.87,141.89) .. controls (526.17,141.89) and (527.22,142.98) .. (527.22,144.33) .. controls (527.22,145.67) and (526.17,146.76) .. (524.87,146.76) .. controls (523.57,146.76) and (522.52,145.67) .. (522.52,144.33) -- cycle ;
\draw  [fill={rgb, 255:red, 0; green, 0; blue, 0 }  ,fill opacity=1 ] (447.49,203.57) .. controls (447.49,202.23) and (448.55,201.14) .. (449.85,201.14) .. controls (451.15,201.14) and (452.2,202.23) .. (452.2,203.57) .. controls (452.2,204.92) and (451.15,206.01) .. (449.85,206.01) .. controls (448.55,206.01) and (447.49,204.92) .. (447.49,203.57) -- cycle ;

\draw (450.85,196.14) node [anchor=north west][inner sep=0.75pt]    {\small $p_{1}$};
\draw (509.46,148) node [anchor=north west][inner sep=0.75pt]    {\small $p_{2}$};

\end{tikzpicture}
		\vspace{-2mm}
		\caption{}
		\vspace{-2mm}
		\label{fig:domainThin}
	\end{subfigure}
	\caption{\small Computational domain: (a) full representation of $\Omega_s$ in $\Omega$, and (b) reduced-dimension problem, 
	 with $\Omega_s$ replaced \marc{by a lower-dimensional region $\Gamma_s$}.}
	\label{fig:domain}
\end{figure}
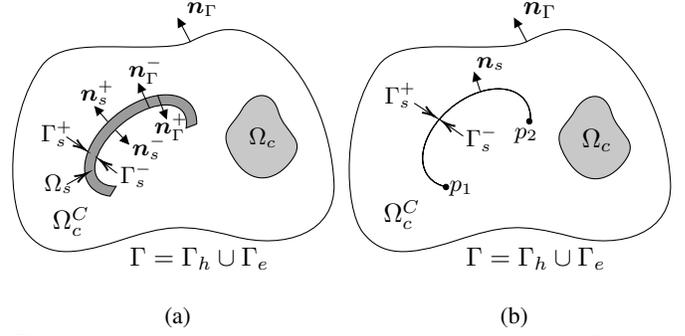

We study the problem of a thin region $\Omega_s$ embedded in a domain $\Omega=\Omega_c \cup \Omega_c^C$, where  $\Omega_c$ and $\Omega_c^C$ denote respectively the conducting and non-conducting parts of $\Omega$. As depicted in Fig.~\ref{fig:domainFull}, the exterior boundary of $\Omega$ ($\partial \Omega = \Gamma$) is composed of two complementary parts $\Gamma_h$ and $\Gamma_e$ (i.e. $\Gamma =\Gamma_h \cup \Gamma_e$ and $\Gamma_h \cap \Gamma_e = \varnothing$) that may be necessary for symmetry or physical purposes such as connecting different subproblems via their common boundaries~\cite{Dular2011}. The thin region $\Omega_s$ belongs to the conductive subdomain ($\Omega_s\subset \Omega_c$) and its interior and exterior boundaries are $\Gamma_s^-$ and $\Gamma_s^+$, respectively. 

\par When \marc{coupling} the TS model with the FEM, $\Omega_s$ is geometrically replaced by a surface located halfway between the original boundaries ($\Omega_{s}$ $\to$ $\Gamma_{s}$ in Fig.~\ref{fig:domainThin}). \marc{In the variational form, by assuming distinct BCs on both sides of $\Gamma_s$, we obtain interface integrals to couple with the TS model within $\Omega_s$.} Thus, the weak form of the $\bm{h}$-formulation, obtained from the weak form of Faraday's law, is defined as follows:
 
\noindent Find $\bm{h} \in \bm{H}(\rm{curl},\Omega)$ such that  
\begin{equation} \label{h_form_weak_}
\begin{split}
&{\Big( {\rho {\rm{ }}\nabla  \times {\bm{h}},\nabla  \times {\bm{h}'}} \Big)_{\textcolor{black}{{\Omega _c}\setminus \Omega_s}}} + {\partial _t}{\Big( {\mu {\rm{ }}{\bm{h}},{\bm{h}'}} \Big)_{\textcolor{black}{\Omega \setminus \Omega_s}} } \\ 
&\resizebox{.88\hsize}{!}{ $ + {\Big\langle {{\bm{n}} \times {\bm{e}},{\bm{h}'}} \Big\rangle _{{\Gamma _{e}}}}
-{\Big\langle {{\bm{n}_s} \times {\bm{e}},{\bm{h}'}} \Big\rangle _{{\Gamma _{s}^+}}}
+ {\Big\langle {{\bm{n_s}} \times {\bm{e}},{\bm{h}'}} \Big\rangle _{{\Gamma _{s}^-}}}
= 0 $},
\end{split}
\end{equation} 
 $\forall$ ${\bm{h}'}$ $\in$  \marc{$\bm{H}_0(\rm{curl},\Omega)$, where ${\bm{h}'}$ are test functions with \mbox{$\bm{n} \times \bm{h}' = 0$ along $\Gamma_h$},} $\bm{n}$ is the outward unit normal vector on $\Gamma$, and ${\left( { \cdot , \cdot } \right)_\Omega }$ and ${\left\langle \cdot , \cdot \right\rangle _\Gamma }$   denote respectively the volume integral over $\Omega$ and the surface integral over $\Gamma$ of the scalar product of their two arguments. Note that \marc{we assume that $\bm{h}\in \bm{H}(\rm{curl},\Omega)$ already satisfies the BC on} ${\bm{n}} \times {\bm{h}}$ along~$\Gamma_h$.

The two last terms of~\eqref{h_form_weak_} express the discontinuity of the tangential components of the electric field along the surface representing the thin region and we considered that \textcolor{black}{\mbox{$\bm{n}_s=-\bm{n}_{\Gamma_s}^+=\bm{n}_{\Gamma_s}^-$} (see Fig.~\ref{fig:domainFull})}. 
Besides these interface terms, the weak form \eqref{h_form_weak_} requires the duplication of the DoFs related to the surface of the thin region. In~\cite{geuzaine2000dual}, the authors propose the decomposition of the field quantities into its continuous and discontinuous parts in order to avoid nodes and edges duplication. This decomposition is also applied in~\cite{sabariego2008h,sabariego2009,gyselinck2008}. 

Here, nodes and edges of the thin surface are duplicated, but except for the nodes located at its extremities (e.g. points $p_1$ and $p_2$ in Fig.~\ref{fig:domainThin}). This creates a {\it crack} in the topological structure, and the non-conducting region becomes non-simply connected. The interfaces $\Gamma_s^+$ and $\Gamma_s^-$ share the nodes at their extremities, such that \mbox{$\Gamma_{s}=\Gamma_{s}^+\cup \Gamma_{s}^-$,} and the tangential components of the magnetic fields on these surfaces are connected by an 1-D FE problem in the thin direction of the sheet.

In order to include the 1-D problem in the weak form~\eqref{h_form_weak_}, we express the surface integral terms on $\Gamma_s$ in~\eqref{h_form_weak_} \marc{by using the variational formulation of Faraday's law over $\Omega_s$, namely}
\begin{equation} \label{volume}
	\begin{split}
-{\Big\langle {{\bm{n}_s} \times {\bm{e}},{\bm{h}'}} \Big\rangle _{{\Gamma _{s}^+}}}
+ {\Big\langle {{\bm{n_s}} \times {\bm{e}},{\bm{h}'}} \Big\rangle _{{\Gamma _{s}^-}}} &= \\ 
{\Big( {\rho {\rm{ }}\nabla  \times {\bm{h}},\nabla  \times {\bm{h}'}} \Big)_{{\Omega _s}}}  + &{\partial _t}{\Big( {\mu {\rm{ }}{\bm{h}},{\bm{h}'}} \Big)_{\Omega_s} }.
	\end{split}
\end{equation}
The volume integrals terms in this expression have opposite signs than those presented in~\eqref{h_form_weak_}, since they are on the right side of~\eqref{volume}. In fact, here we should consider the outward unit normal vector of \textcolor{black}{the boundary of $\Omega_{s}$,} i.e. $\bm{n}_s^\pm$ \textcolor{black}{in Fig.~\ref{fig:domainFull}, but for the sake of simplicity, we denoted $\bm{n}_s^+=-\bm{n}_s^-=\bm{n}_s$(Fig.~\ref{fig:domainThin}).}

Inside the sheet, we assume that the local magnetic field is written as ${\bm{h}}_x (x,y,z,t)={\bm{h}}_x(x,z,t) \zeta(y)$ and the test function as ${\bm{h}}_x' (x,y,z)={\bm{h}}_x'(x,z) \zeta'(y)$, with ${\bm{h}}_x (x,z,t)$ and ${\bm{h}}_x' (x,z,t)$ tangential to $\Gamma_s$, and $\zeta(y)$ and $\zeta'(y)$ differentiable \textcolor{black}{in the interval} $-d/2 \leq y\leq d/2$. The volume integrals terms in~\eqref{volume} are then reduced to surface integrals terms as follows
\vspace{-1mm}
\begin{equation} \label{term1}
\begin{split}
&\resizebox{.89\hsize}{!}{ ${\Big( {\rho {\rm{ }}\nabla  \times {\bm{h}},\nabla  \times {\bm{h}}'} \Big)_{{\Omega _s}}} = \Big( \rho \nabla \times ({\bm{h}}_x \zeta) , \nabla \times ({\bm{h}}_x' \zeta') \Big)_{\Omega_s} $} \\
&\resizebox{.89\hsize}{!}{$= \Big( \rho \left( \textcolor{black}{\zeta} \nabla \times {\bm{h}}_x -  {\bm{h}}_x \times \nabla \zeta \right), \zeta' \nabla \times {\bm{h}}_x' - {\bm{h}}_x' \times \nabla \zeta' \Big)_{\Omega_s} $}\\
&\myeq \Big\langle  {\bm{h}}_x,{\bm{h}}_x'  \Big\rangle_{\Gamma_s} \cdot \int_{-d/2}^{d/2} \rho \partial_y \zeta \partial_y \zeta' dy,
\end{split}
\end{equation}

\begin{equation} \label{term2}
\begin{split}
{\partial _t}{\Big( {\mu {\rm{ }}{\bm{h}},{\bm{h}}'} \Big)_{\Omega_s} } &= \partial_t \Big(\mu {\bm{h}}_x \zeta, {\bm{h}}_x' \zeta' \Big)_{\Omega_s} \\
 &= \partial_t \Big\langle  {\bm{h}}_x, {\bm{h}}_x' \Big\rangle_{\Gamma_s} \cdot \int_{-d/2}^{d/2} \mu \zeta \zeta' dy,
\end{split}
\vspace{-2mm}
\end{equation}
where $\Gamma_s = \Gamma_{s}^+ \cup \Gamma_{s}^-$. Note that the expansion in~\eqref{term1} has been reduced to the 2-D case, so that ${\bm{h}}_x$ is independent of~$z$, i.e., ${\bm{h}}_x(x,t)$, and the terms $\nabla\times{\bm{h}}_x$ and $\nabla\times{\bm{h}}'_x$ vanish. In 3-D, these terms should be taken into account. 

By choosing $\zeta$ and $\zeta'$ as $\theta_p$ and $\theta_q$, respectively, the integral terms in~\eqref{term1} and~\eqref{term2} \marc{are seen to be components of the elementary matrices} $[\mathcal{S}]$ and $[\mathcal{M}]$ in~\eqref{Spq} and~\eqref{Mpq}. Finally, \marc{taking} the weak form of the \mbox{$\bm{h}$-formulation}~\eqref{h_form_weak_} with the interface terms rewritten~\eqref{volume} using~\eqref{term1} and~\eqref{term2}, \marc{then estimating for each degree of freedom in $\Gamma_{s}(x,z)$ the variation in $y$ using~\eqref{H_inside} with a system of the form~\eqref{EqSys},  we obtain a coupled system of equations for the magnetic field inside and outside the TS}. \textcolor{black}{The IBCs in the proposed TS model are obtained from~\eqref{volume}-\eqref{term2}. Note that, with a single pair of hyperbolic basis functions, these equations become equivalent to the IBCs~\eqref{IBC01} and~\eqref{IBC02} of the classical TS model.}

\section{Validation and Application} \label{Validation}

\begin{figure}[b]
	\vspace{-5mm}
	\centering
	\import{}{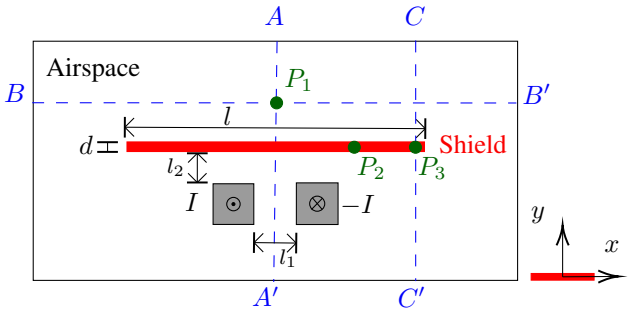}
	\caption{\small Geometry of the planar shield placed over a pair of wires, and lines $AA'$ $BB'$ \textcolor{black}{and $CC'$}, and points $P_1$, $P_2$ \textcolor{black}{and $P_3$} where the local distributions of the fields are analyzed.} 
	\label{fig:shield}
\end{figure}

We consider a 2-D planar shield (width \mbox{$l= 1$\,m} and thickness \mbox{$d=1$\,mm}) placed over a pair of wires carrying a current $\pm I$ (Fig.~\ref{fig:shield}). The conductors are 2x2\,cm$^2$ separated by a distance of $l_1=30$\,cm, and the distance between the conductors and the shield is $l_2=10$\,cm. The free-space region is 4x4\,m$^2$. The coordinate system $xyz$ is defined at the center of the shield geometry and we evaluate the magnetic field distribution along the lines \mbox{$AA'$($x=0$, $y$),} \mbox{$BB'$($x,y=10$\,cm)} \textcolor{black}{and \mbox{$CC'$($x=l/2-l/100,y$),}} and at points \mbox{$P_1$($x=0$,$y=10$\,cm),} \mbox{$P_2$($x=l/4$,0)} \textcolor{black}{and \mbox{$P_3$($x=l/2-l/100$,0)}}. 

The application of standard FEM using the $\bm{h}$-formulation with a full 2-D representation of the shield gives the reference solution. The solutions obtained with the application of the TS model are here compared with the reference solution in terms of local field distributions, Joule losses and mesh simplification.

\textcolor{black}{The relative difference ($\mathcal{R}$) between the solutions is calculated as}
\begin{equation}\label{error}
\textcolor{black}{\mathcal{R}\text{ } [\%]= \frac{\| \text{TSS}- \text{FES} \|_2 }{ \| \text{FES} \|_2 }\times 100,}
\end{equation}
\textcolor{black}{where $\text{TSS}$ and $\text{FES}$ are the TS and the reference FE solutions, respectively, and} $\|.\|_2$ denote the Euclidean norm of the argument.

\begin{figure}[t]
	\begin{subfigure}{0.49\textwidth}
 		\fbox{\includegraphics[trim={0cm 5cm 15cm 5cm}, clip, width=0.48\textwidth]{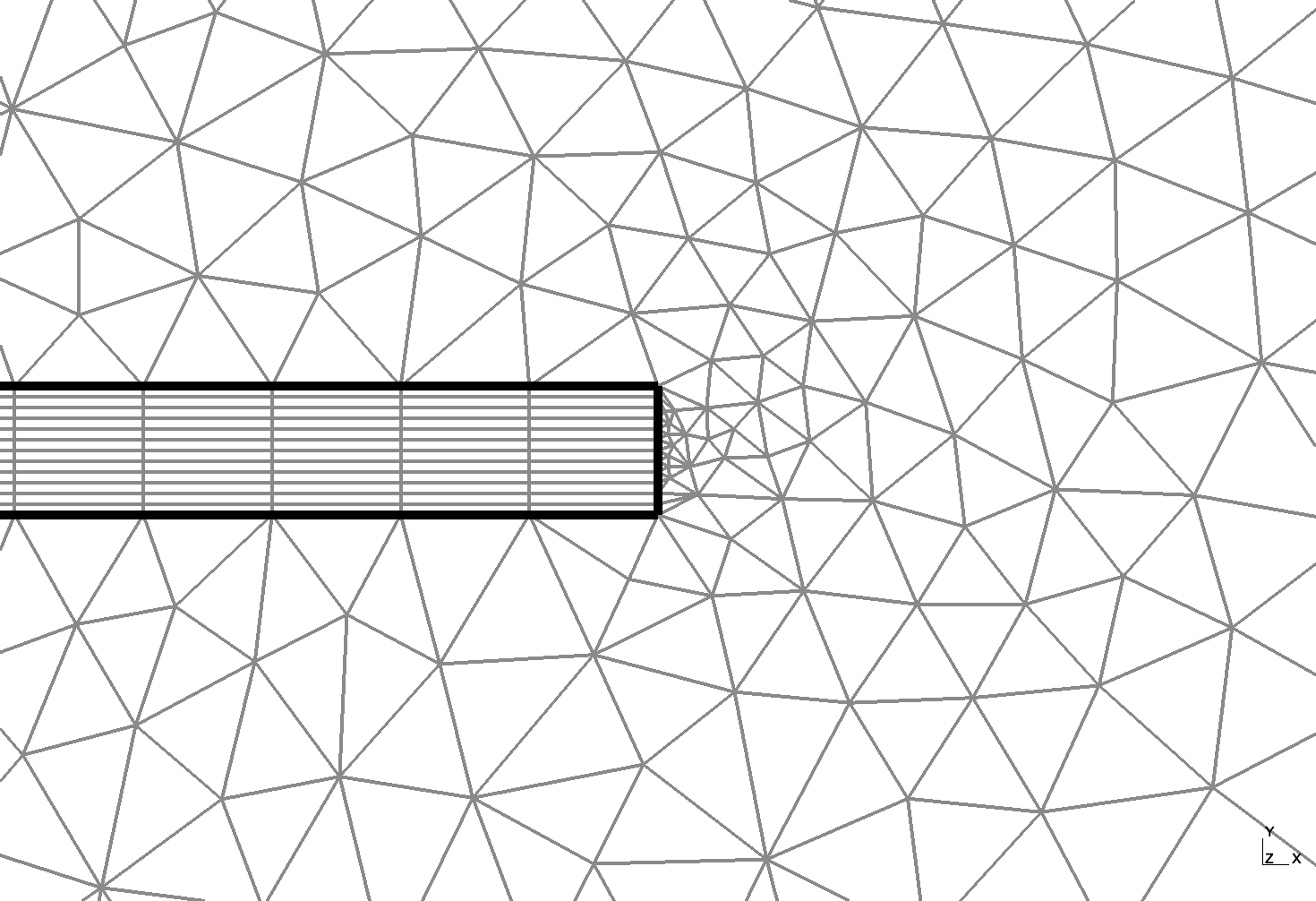}}
 		\fbox{\includegraphics[trim={0cm 5cm 15cm 5cm}, clip, width=0.48\textwidth]{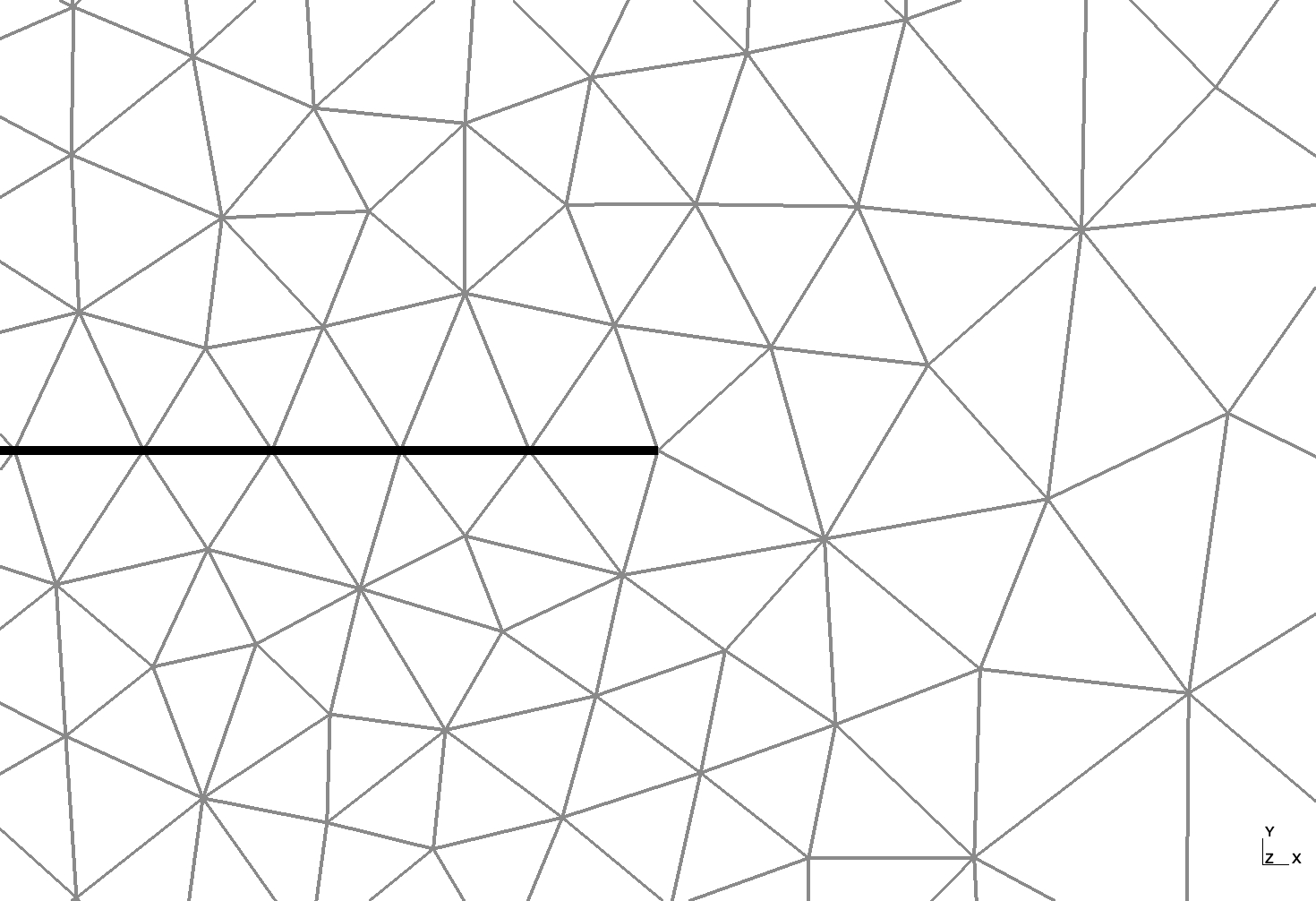}}
 		\caption{Mesh with 1~mm wide elements: FE (left) and TS (right).}
 		\label{mesh1}
 	\end{subfigure}
	\begin{subfigure}{0.49\textwidth}
		\fbox{\includegraphics[trim={0cm 5cm 15cm 5cm}, clip, width=0.48\textwidth]{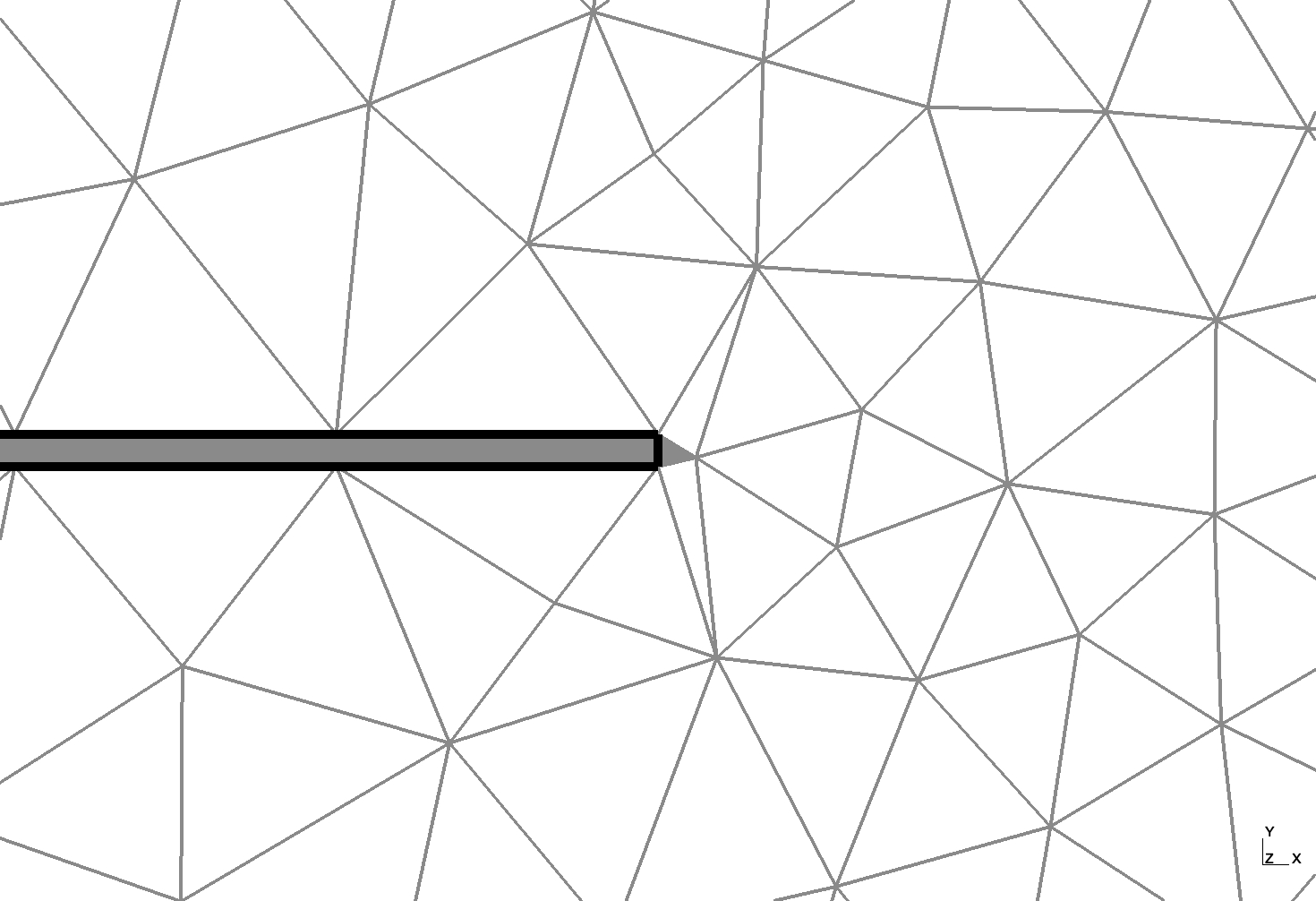}}
		\fbox{\includegraphics[trim={0cm 5cm 15cm 5cm}, clip, width=0.48\textwidth]{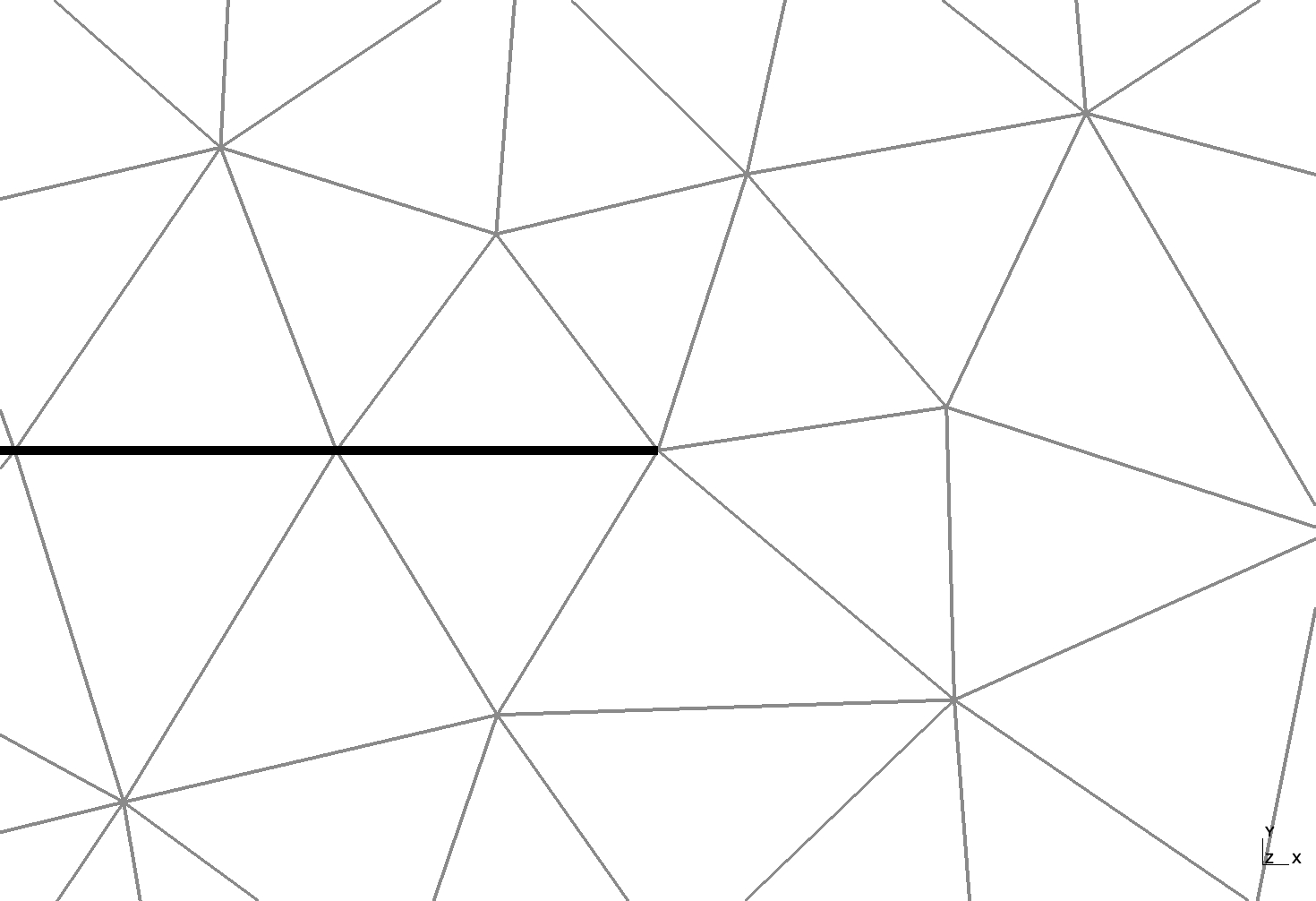}}
		\caption{Mesh with 10~mm wide elements: FE (left) and TS (right).}
		\label{mesh2}
	\end{subfigure} 
	\caption{\small \textcolor{black}{Mesh differences near the right edge of the shield with (a) 1\,mm and (b) 10\,mm wide elements. For better visualization of the elements in the surroundings of the shield, figures in (a) and (b) are not to same scale. Figures in (a) were zoomed-in $\times4$ compared to figures in (b).}}
	\label{mesh}
	\vspace{-5mm}
\end{figure}

In terms of mesh parameters, we defined a structured rectangular mesh in $\Omega_s$ with 12 elements across the shield thickness in the FE model \textcolor{black}{(Fig.~\ref{mesh1}-left)}. The shield surface was discretized in 1\,mm wide elements, and 100\,mm wide elements were considered on the external boundary. Moreover, first-order \textcolor{black}{edge elements} were used \textcolor{black}{in $\Omega$}.

With the described mesh configuration, the application of the TS model \textcolor{black}{(Fig.~\ref{mesh1}-right})  with $n=1$ represents a reduction in the total number of DoFs by 14.3\% in comparison with the FE model. However, a more significant reduction can be achieved with the TS model since a coarser mesh can be employed while maintaining a sufficiently high mesh quality \textcolor{black}{and solution accuracy.}

According to~\cite{marchandise2014optimal}, the quality of the triangular mesh can be evaluated by computing the aspect ratio of the inscribed radius to the circumscribed radius of every triangle. \textcolor{black}{For the meshes presented in Fig.~\ref{mesh1},} the smallest aspect ratios are 0.408 and 0.513 for the 2-D FE and the TS models, respectively. If a coarser mesh is considered, e.g., elements of size $10$\,mm in the shield surface \textcolor{black}{(Fig.~\ref{mesh2}),} these aspect ratios become 0.085 and 0.621, respectively. In this case, the low quality of the mesh in the standard FE may reduce the solution accuracy and its convergence. Therefore, with elements of 10~mm, the TS approach is preferable. It allows reducing the number of DoFs by 80.9\% \textcolor{black}{compared to mesh in the FE reference model} while maintaining the initial mesh quality in the 2-D domain. Thus, 10~mm wide elements were \textcolor{black}{used} in the TS model.

All the models described in this paper were implemented in the open-source code Gmsh~\cite{Gmsh2009} and the solver GetDP~\cite{GetDP2013}. Simulations were conducted on a personal computer with an Intel i7 2400 processor and 16 Gb of memory.  Harmonic and time-transient simulations for different types of shields were performed.

\subsection{Time-Harmonic Regime}

Simulations were first performed in the harmonic regime for the sake of validation of the proposed TS \textcolor{black}{model}. The current~$I$ was set to \textcolor{black}{6\,kA} (current density of 15\,A/mm$^2$) at an operating frequency of~$f= 50$\,Hz, and two shield configurations were studied:

\noindent {{\bfseries{shield 1}}: $\mu_r=1$ and $\sigma=1$\,MS/m $\Rightarrow$ $\delta= 71.2$\,mm,}

\noindent {{\bfseries{shield 2}}: $\mu_r=1000$ and $\sigma=10$\,MS/m $\Rightarrow$ $\delta = 0.712$\,mm,}

\noindent where $\mu_r=\mu/\mu_0$ is the relative magnetic permeability with respect to $\mu_0$, the magnetic permeability of air.

The hyperbolic basis functions $\psi^\pm$ in~\eqref{psi_basis} were defined accordingly. With the first shield configuration, the basis functions are equivalent to those presented in Fig.~\eqref{Basis_deltaggd}, since $\delta \gg d$. With the second configuration, the hyperbolic basis functions were defined with $\delta = 0.712$\,mm (or $\delta=1.40 d$, since $d=1$\,mm). For harmonic regime simulations, only one pair of \textcolor{black}{hyperbolic} basis functions was considered ($n=1$).

\begin{figure} [t]
	\hspace{-0.3cm}
	\includegraphics[trim={0cm 0cm 0cm 0cm},clip,width=0.5\textwidth]{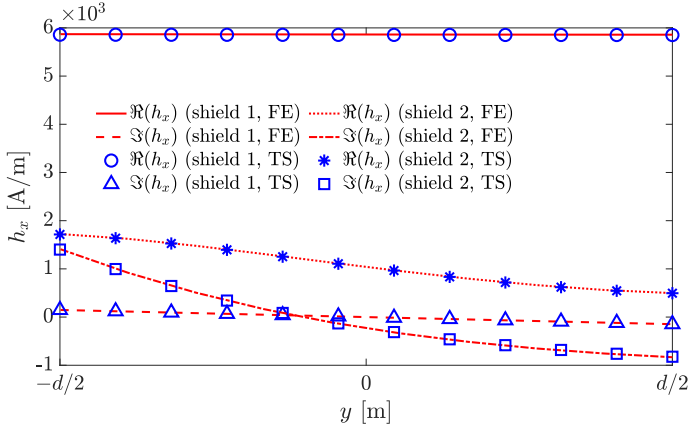}
	\caption{\small Profile of~$h_x$ at point~$P_2$ inside shield~1~($\delta \gg d$) and shield~2~($\delta < d$). \textcolor{white}{AA} }
	\vspace{-0.5cm}
	\label{FieldInHarmonic}
\end{figure}

\begin{figure} [t]
	\hspace{-0.3cm}
	\includegraphics[trim={0cm 0cm 0cm 0cm},clip,width=0.5\textwidth]{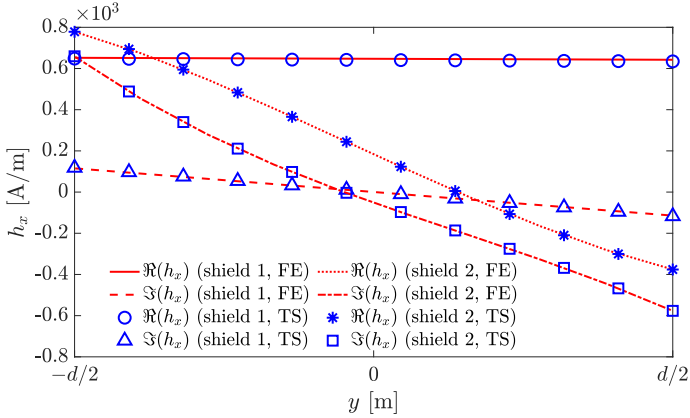}
	\caption{\small Profile of~$h_x$ at point~$P_3$ inside shield~1~($\delta \gg d$) and shield~2~($\delta < d$). \textcolor{white}{AA} }
	\vspace{-0.5cm}
	\label{FieldInHarmonic100}
\end{figure}

Real and imaginary components of~$h_x$ at points~$P_2$ \textcolor{black}{and~$P_3$} inside the shields are presented in Fig.~\ref{FieldInHarmonic} \textcolor{black}{and in Fig.~\ref{FieldInHarmonic100}}, respectively. The field profiles from the TS model were obtained by evaluating the field intensities in the shield with the proposed hyperbolic basis functions. \textcolor{black}{Excellent} agreement with reference solutions in terms of real and imaginary components of $h_x$ were observed for both shield configurations. \textcolor{black}{The excellent agreement at point $P_3$ shows that the proposed model can provide accurate solutions also near the extremities of the shield. Therefore, even though no special consideration has been made at its extreme points, edge effects are correctly represented.}
 
 \begin{figure}[t]
 	\vspace{3mm}
 	\begin{subfigure}{0.49\textwidth}
 		\includegraphics[trim={0cm 10cm 15cm 0cm}, clip, width=0.48\textwidth]{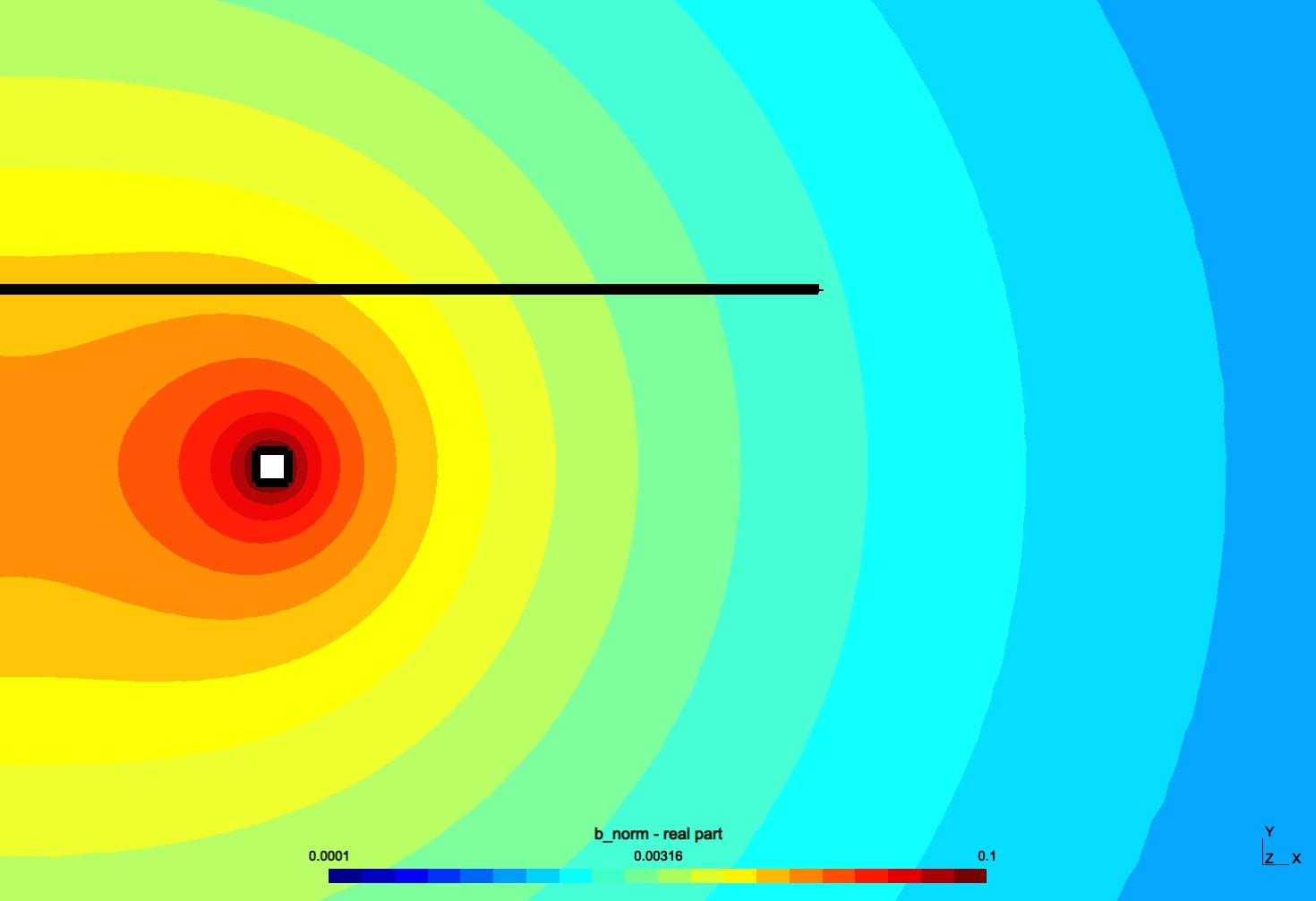}
 		\includegraphics[trim={0cm 10cm 15cm 0cm}, clip, width=0.48\textwidth]{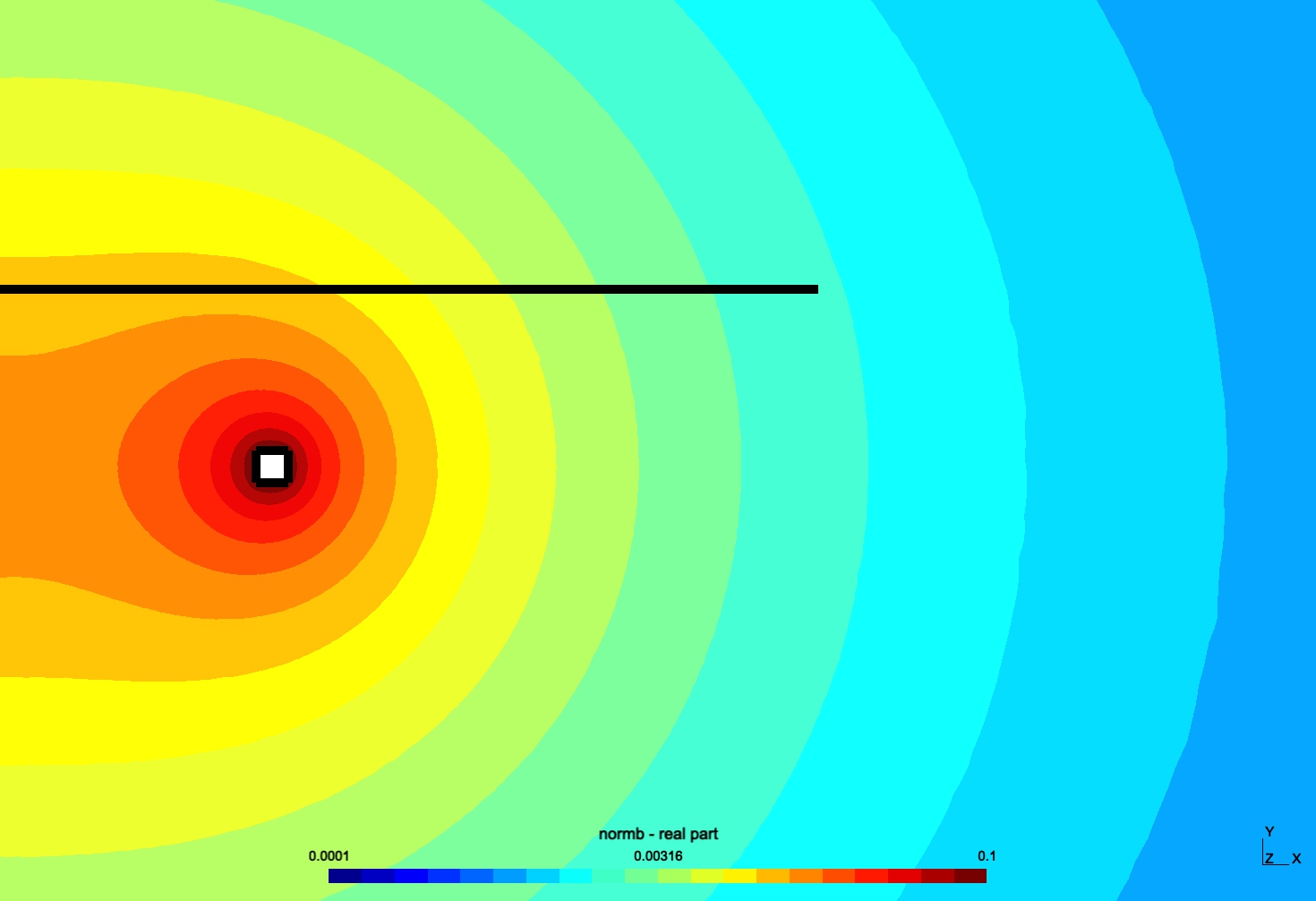}
 		\includegraphics[trim={12cm 0.5cm 12cm 33.3cm}, clip, width=1\textwidth]{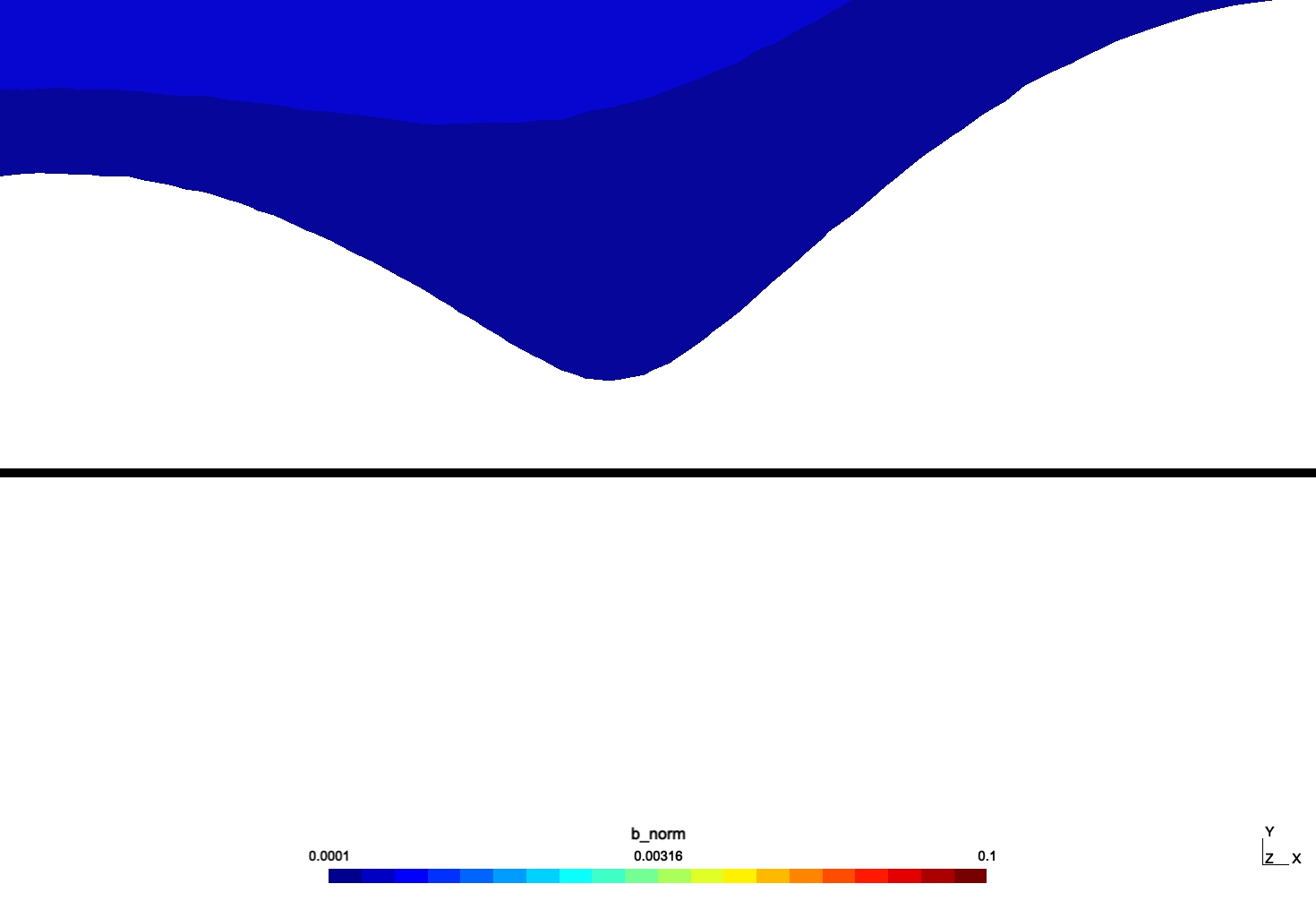}
 		\caption{shield 1: FE solution (left) and TS solution (right) \vspace{2mm}}
 		\label{fluxshield1}
 	\end{subfigure} 
 	\begin{subfigure}{0.49\textwidth}
 		\includegraphics[trim={0cm 10cm 15cm 0cm}, clip, width=0.48\textwidth]{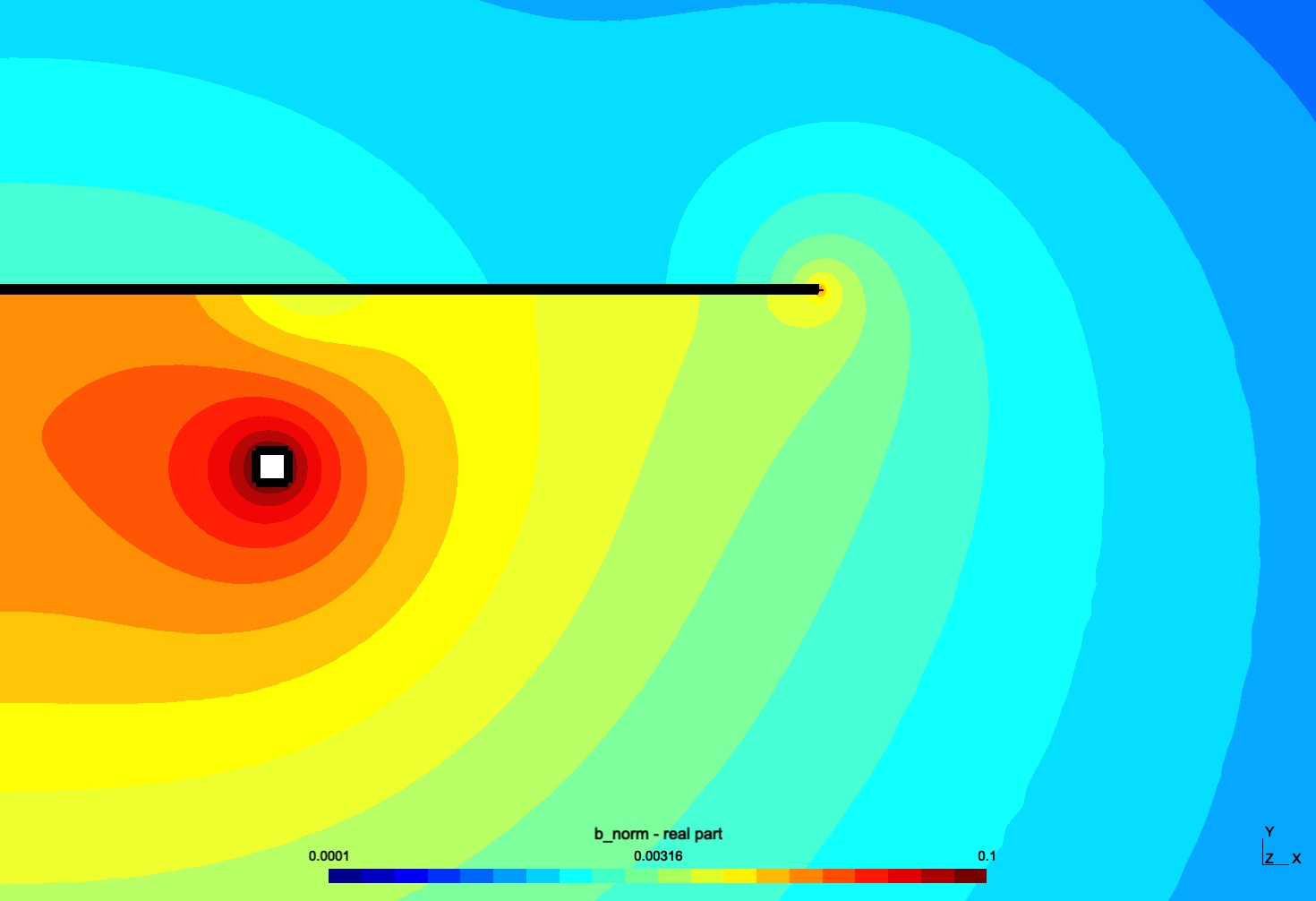}
 		\includegraphics[trim={0cm 10cm 15cm 0cm}, clip, width=0.48\textwidth]{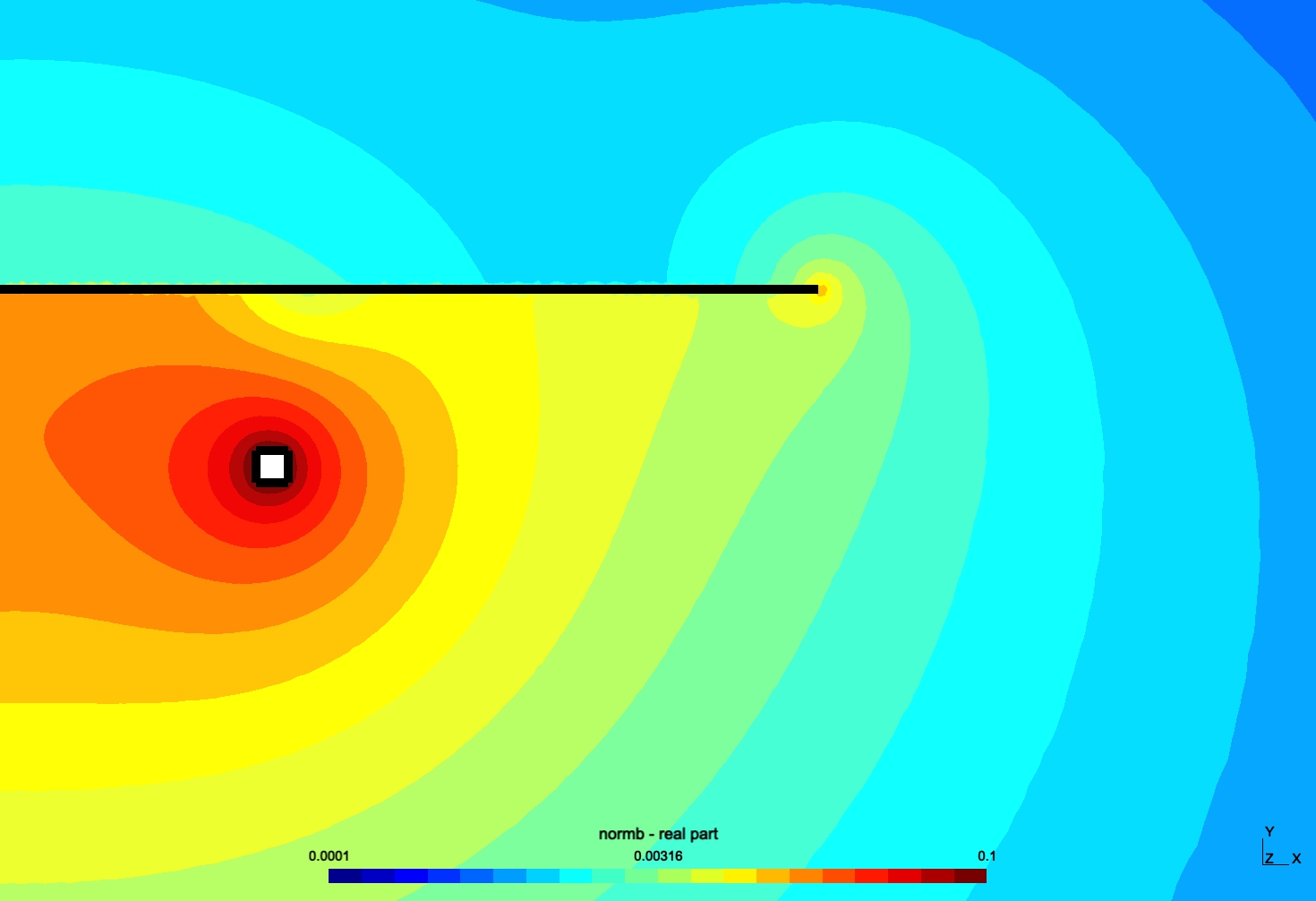}
 		\includegraphics[trim={12cm 0cm 12cm 33.3cm}, clip, width=1\textwidth]{FiguresA/Harmonic/Legend.jpg}
 		\caption{shield 2: FE solution (left) and TS solution (right)}
 		\label{fluxshield2}
 	\end{subfigure}
 	\caption{\textcolor{black}{Shaded plot of the isovalues of the magnetic flux density ($|\bm{b}|$) in half of domain (a) shield 1, and (b) shield 2. Note the difference between the solutions depending on the shielding configuration, and the equivalence between the FE (right) and the TS (left) solutions in both cases. The airspace is not to scale.}}
 	\label{fluxshield}
 \end{figure}

\textcolor{black}{In Fig.~\ref{fluxshield}, we present a shaded plot of the magnetic flux density.  With the first shield configuration (Fig.~\ref{fluxshield1}), the replacement of the original 2-D region by a thin sheet has no noticeable impact on the magnetic flux density distribution.  
	Indeed, from Fig.~\ref{FieldInHarmonic} and Fig.~\ref{FieldInHarmonic100}, we observe that $h_x$ is almost constant inside shield~1. The tangential components of the field are continuous on the surface representation in the TS model. However, with the second shield configuration (Fig.~\ref{fluxshield2}), the physics inside the plate produces a discontinuity of $h_x$ that deforms the flux lines in the air surrounding the edge. Thus, in both shield configurations, the solution from the TS model agrees with the FE solution in terms of field distributions inside and outside the shield. }

\begin{table}[t]
		\caption{\small \textcolor{black}{Relative differences for the magnitude of the magnetic field on lines $AA'$, $BB'$, $CC'$, and points $P_2$ and $P_3$ depicted in Fig.~\ref{fig:shield}.} }
	\centering
	\begin{tabular}{@{}llllll@{}}
		\toprule
		& $AA'$  & $BB'$  & $CC'$  & $P_2$  & $P_3$  \\ \midrule
		Shield 1 & 0.72\% & 1.22\% & 1.95\% & 0.14\% & 0.92\% \\
		Shield 2 & 0.90\% & 1.56\% & 2.46\% & 1.29\% & 2.46\% \\ \bottomrule
	\end{tabular}
\vspace{-0.5cm}
\label{Rdiff}
\end{table}

The relative differences for the local magnetic field along lines $AA'$, $BB'$ \textcolor{black}{and $CC'$}, \textcolor{black}{and at points $P_2$ and $P_3$} were computed with~\eqref{error}. \textcolor{black}{The \mbox{$\mathcal{R}$-values} are summarized in Table~\ref{Rdiff}. The maximum difference is 2.46\% and occurs on line $CC'$ and point $P_3$ with shield 2. This difference may be related to the geometrical difference between the TS and the FE models. Since the thickness of the shield is not represented in the TS model, it is expected to observe at least a slight difference near its extremities.}

\subsection{Time-Transient Regime (Linear Case)}

\begin{figure} [t]
	\hspace{-0.3cm}
	\includegraphics[trim={0cm 0cm 0cm 0cm},width=0.5\textwidth]{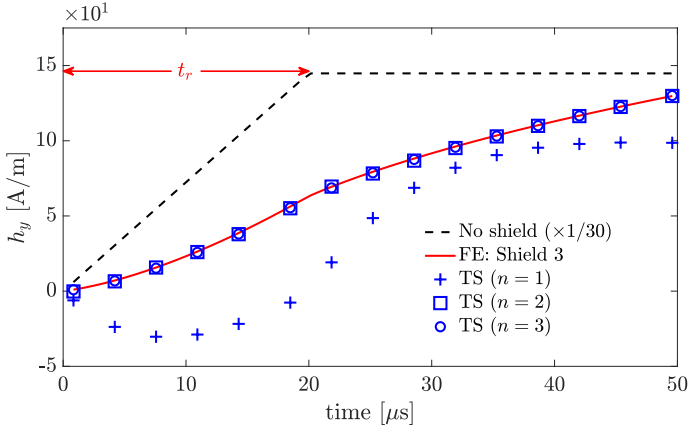}
	\caption{\small Time-evolution of $h_y$ at point $P_1$ with pulsed current imposed to the wires and a ferromagnetic shield configuration (Shield~3). Results obtained with $n$ up to 3 in the TS model compared with the FE solution.}
	\label{fig:pulsedFerro}
\end{figure}

\begin{figure} [t]
	\hspace{-0.3cm}
	\includegraphics[trim={0cm 0cm 0cm 0cm},width=0.501\textwidth]{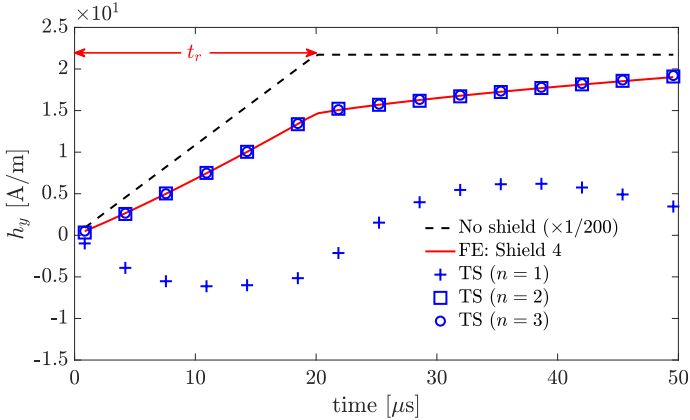}
	\caption{\small Time-evolution of $h_y$ at point $P_1$ with pulsed current imposed to the wires and a conductive shield configuration (Shield~4). Results obtained with $n$ up to 3 in the TS model compared with the FE solution.}
	\label{fig:pulsedCond}
	\vspace{-3mm}
\end{figure}

Time-transient simulations of the planar shield were also performed. In this study, we applied a pulsed current source in the wires whose waveform produces the magnetic field shown in~Fig.~\ref{fig:pulsedFerro} and~\ref{fig:pulsedCond} (black dashed lines) in absence of the shield. The amplitude of the current~$I$ is the same as in the harmonic regime cases (\textcolor{black}{$|I|$= 6\,kA}), and the rise time ($t_r$) was set \mbox{to 20\,$\mu$s.} The simulation time was $t_{\text{ max}}= 50$\,$\mu$s with a time-step of $t_{\text{ max}}/120$. The implicit Euler scheme was used.

In terms of material composition, two new shielding configurations were studied:

\noindent {{\bfseries{shield 3}}: $\mu_r=1000$ and $\sigma=1$\,MS/m,}

\noindent {{\bfseries{shield 4}}: $\mu_r=100$ and $\sigma=10$\,MS/m.}

Shields 3 and 4 are both ferromagnetic, but shield 4 is more conductive than shield 3. \textcolor{black}{The simulation time was chosen to define the fundamental frequency as $f=1/(4t_\text{max})$, i.e. $f = 5$\,kHz. Therefore, these shields have same penetration depth~$\delta = d/4.44$.} Consequently, the same hyperbolic basis functions can be used to tackle these problems. The first set of basis functions was defined with frequency equal to the fundamental ($f_1=f$). Additional basis functions were then calculated using odd harmonic frequencies of $f_1$, i.e., $f_k/f_1=2k-1$, with $1\le k \le n$. 

In Fig.~\ref{fig:pulsedFerro} and~\ref{fig:pulsedCond}, the time-evolution of $h_y$ at point $P_1$ was compared with the reference solution for shield 3 and 4, respectively. The number of basis functions in the TS model was varied from 1 to 3, and the solutions approached the reference solution as $n$ increased. 

\begin{figure} [t]
	\vspace{3mm}
	\hspace{-0.3cm}
	\includegraphics[trim={0cm 0cm 0cm 0cm},width=0.5\textwidth]{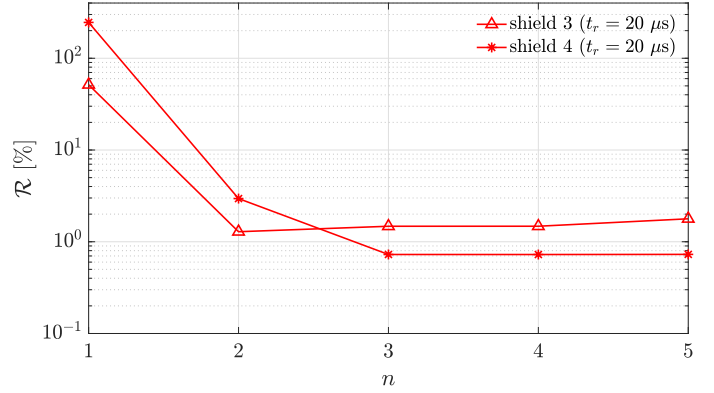}
	\caption{\small Relative difference $\mathcal{R}$ of instantaneous $h_y$ at point $P_1$ as a function of $n$ in the TS approach for \textcolor{black}{ $t_r=20$\,$\mu$s} in the time-transient study. \textcolor{black}{Note the fast convergence towards the FE solution.}}
	\label{fig:pulsed_error}
\end{figure}

The relative differences of the profiles of $h_y$ in the TS model to the reference solution as a function of the number of basis functions $n$ are presented in Fig.~\ref{fig:pulsed_error}. The maximum relative difference~$\mathcal{R}$ decreases from 245.8\% with $n=1$ to 2.95\% with $n=2$, and to less than 2\% for $n\geq3$. Simulations with $n> 3$ show little improvement in terms of solution accuracy. \textcolor{black}{This is mainly due to the geometrical differences between the FE and the TS models, as discussed in the time-harmonic regime case. Despite this, the application of the TS model shows a good compromise between computational cost and solution accuracy.}

The number of DoFs, the computation time and the Joule losses in the reference and TS models are summarized in Table~\ref{Table1}. Since a coarser mesh was applied in the TS model, simulations are more than five times faster with this approach than with standard FE. \textcolor{black}{The number of DoFs is nearly independent of~$n$,} and is reduced by more than 70\%, even with $n=5$. Despite this, the relative error in losses estimation is less than 2\% with $n\geq3$ in shields 3 and 4.

\begin{table}[t]
	\caption{\small Number of DoFs, CPU time and total Joule losses in shield 3 and 4 in time-transient analysis. }
	\centering
	\begin{tabular}{@{}lcccc@{}}
		\toprule
		\multicolumn{1}{c}{Model} & \begin{tabular}[c]{@{}c@{}}Number \\ of DoFs\end{tabular} & \begin{tabular}[c]{@{}c@{}}CPU\\ time {[}s{]}\end{tabular} & \begin{tabular}[c]{@{}c@{}}Joule losses\\ Shield 3 \\ {[}J/m{]}\end{tabular} & \begin{tabular}[c]{@{}c@{}}Joule losses\\ Shield 4\\ {[}J/m{]}\end{tabular} \\ \midrule
		Standard FE                    & 176054                                                    & 1040.12                                                    & 2.4915                                                                   & 0.7514                                                                   \\
		TS ($n=1$)                & 33636                                                     & 189.97                                                     & 2.4503                                                                    & 0.7676                                                                   \\
		TS ($n=2$)                & 34036                                                     & 205.19                                                     & 2.5103                                                                    & 0.7596                                                                   \\
		TS ($n=3$)                & 34436                                                     & 221.40                                                     & 2.5099                                                                    & 0.7592                           \\                                        
		TS ($n=4$)                & 34836                                                     & 276.58                                                     & 2.5099                                                                    & 0.7592                                                                   \\
		TS ($n=5$)                & 35236                                                     & 309.15                                                     & 2.5087                                                                    & 0.7592                                                                   
		\\ \bottomrule
	\end{tabular}
	\label{Table1}
\end{table}

\vspace{-3mm}
\subsection{Time-Transient Regime (Nonlinear Case)}

\textcolor{black}{In an attempt to extend the proposed TS model to nonlinear analysis,} the shielding problem was also studied \marc{for shields with} nonlinear magnetic permeability ($\mu = \mu(\bm{h})$). A sinusoidal supply current of amplitude \textcolor{black}{$|I|=6$\,kA} at an operating frequency $f=1$\,kHz was applied to the wires. The effects of the saturation and the influence of the number basis functions in the proposed TS model were analyzed in terms of solution accuracy.

\begin{figure} [t]
	\vspace{0.3cm}
	\hspace{-0.3cm}
	\includegraphics[trim={0cm 0cm 0cm 0cm},width=0.5\textwidth]{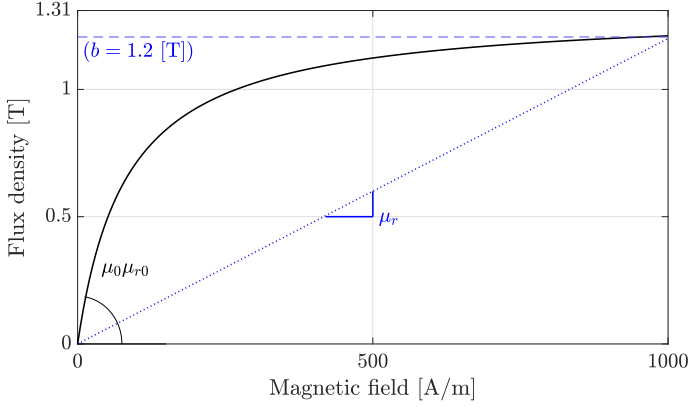}
	\caption{\small \mbox{$B$-$H$} \textcolor{black}{saturation} curve obtained from~\eqref{BHrelation} with \mbox{$ \mu_0 m_0=1.31$} and \mbox{$\mu_{r0} = $ 12500.} The intersection of the \mbox{$B$-$H$} curve  with the horizontal dashed line gives $\mu_r=$ 1000,  which was used \marc{to parametrize} the hyperbolic basis functions. }
	\label{BHcurve}
	\vspace{-3mm}
\end{figure}

\marc{The material properties were modelled with} an isotropic saturation law \marc{expressing} the magnetic permeability as a function of the magnetic field intensity, i.e.,
\vspace{-1mm}
\begin{equation} \label{BHrelation}
\mu (\bm{h})= \mu_0\left( 1+ \left( \frac{1}{\mu_{r0}-1}+ \frac{|| \bm{h}||}{m_0}\right)^{-1} \right),
\vspace{-1mm}
\end{equation}
where $\mu_{r0}$ is the relative permeability at origin and $m_0$ the saturation magnetic field in A/m. The differential permeability required for the application of the NR-scheme was defined as in~\cite{Dular2020}.
We carried out simulations with $\mu_0 m_0=1.31$ and \mbox{$\mu_{r0}=12500$}.
The $B$-$H$ curve is presented in Fig.~\ref{BHcurve}. Furthermore, the electrical \textcolor{black}{conductivity} of the shield was fixed at \mbox{$\sigma=1$\,MS/m.}

One time period was simulated, i.e., $t_{max}=1/T$, with \mbox{$T=1/f$.} Moreover, the time-step was set to $\Delta t = T/120$, and the maximum number of iterations for the NR-scheme was set to 12 in both the reference and the proposed TS models. Besides, the number of points used in the Legendre-Gauss quadrature was 20 points. This number of points is considered sufficiently high to avoid errors related to the numerical integration of the hyperbolic functions across \textcolor{black}{whole} the thickness of the shield. \textcolor{black}{ Depending on the penetration depths of the basis functions, less integration points could be used, but we kept 20 at all times to remains on the safe side.}

\begin{figure} [t]
	\hspace{-0.3cm}
	\includegraphics[trim={0cm 0cm 0cm 0cm},width=0.5\textwidth]{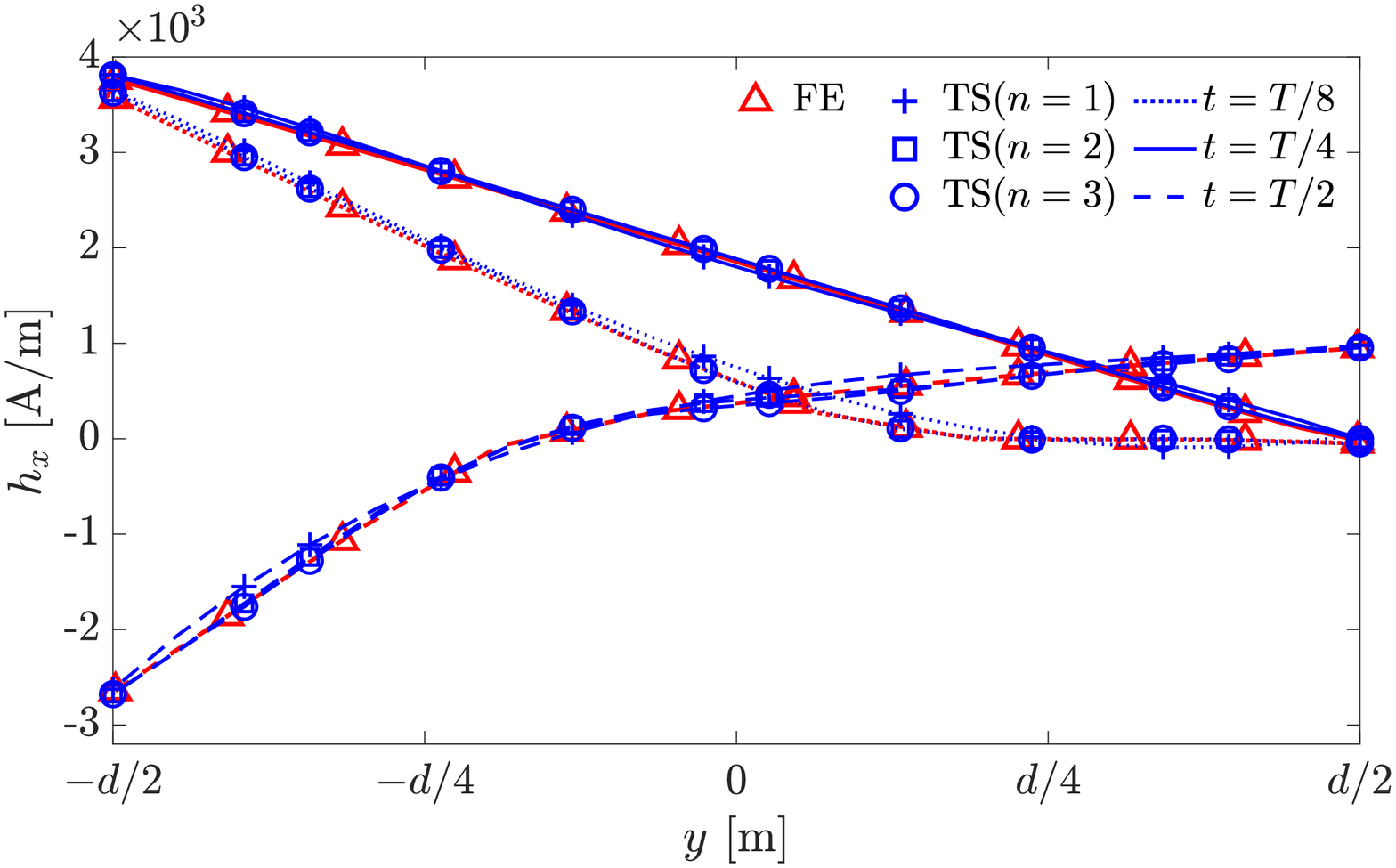}
	\label{MagFieldIn3}
	\vspace{-0.2cm}
	\begin{tikzpicture}[overlay]
	\node[] at (4.42,2.22) {\includegraphics[trim={0cm 0cm 0cm 0cm},width=0.145\textwidth]{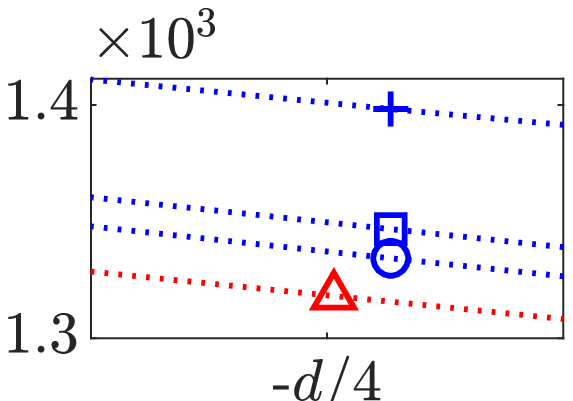}};
	\end{tikzpicture}
	\begin{tikzpicture}[overlay]
	\node[] at (7.0,2.13) {\includegraphics[trim={0cm 0cm 0cm 0cm},width=0.145\textwidth]{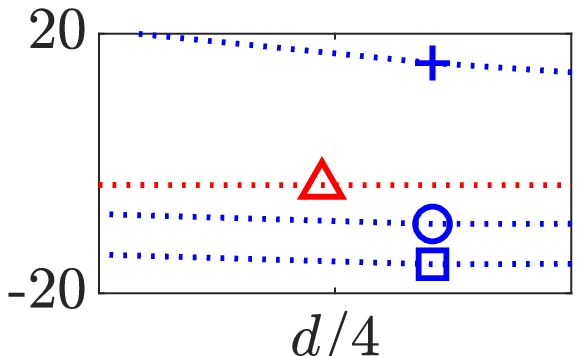}};
	\end{tikzpicture}
	\vspace{-0.3cm}
	\caption{\small Profiles of $h_x$ inside the nonlinear shield at point $P_2$ and at $t=T/8$, $T/4$ and $T/2$ obtained with the 2-D FE solution and with the TS model for $n=1$ to 3. The insets show the solutions at $T/8$ and $y=\pm d/4$.}
	\label{InteriorFieldNolinear}
	\vspace{-0.2cm}
\end{figure} 

\begin{figure} [t]
\textbf{}	\includegraphics[trim={0cm 0cm 0cm 0cm},width=0.49\textwidth]{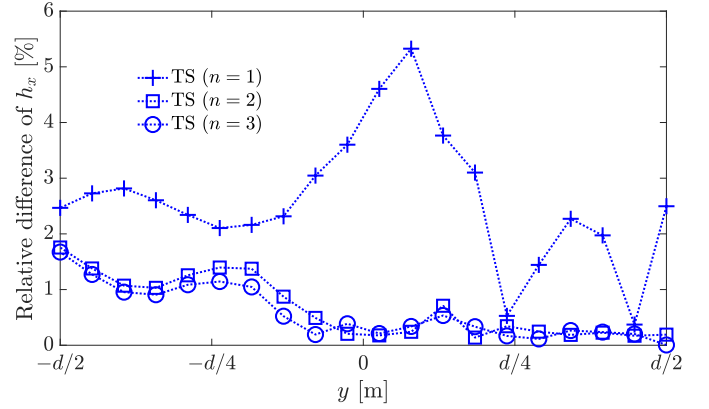}
	\caption{\small \textcolor{black}{Relative difference between the TS and the reference FE solutions for the $h_x$ profile inside the nonlinear shield at $P_2$ and $t=T/8$.}}
	\label{InteriorFieldNolinearError}
	\vspace{-3mm}
\end{figure}

The first set of hyperbolic basis functions was defined by taking  $f_1=f$ and higher order basis functions that are odd multiples of $f_1$. Furthermore, the magnetic permeability used in the definition of the basis functions was taken from the $B$-$H$ curves corresponding to a flux density \mbox{$b=1.2$ T,} i.e., \mbox{$\mu_r=$ 1000} for \mbox{$\mu_{r0}=12500$}. A similar approach was used in~\cite{Sabariego2010} \marc{to parametrize} basis functions for nonlinear SIBCs.

Fig.~\ref{InteriorFieldNolinear} shows the $h_x$ profile throughout the thickness of the shield for three specific simulation times ($T/8$, $T/4$ and $T/2$). 
Results are presented for $n=1$ to 3 and compared with the 2-D FE solution. Note that the accuracy of the proposed TS model clearly improves with $n$. \textcolor{black}{Since the penetration depth used to define the hyperbolic basis functions is at the same time inversely proportional to the square root of $f$ and $\mu_r$, the additional frequency components can be interpreted as a way to consider an increase of the magnetic permeability. For this reason, the saturation effects observed at $t=T/8$ and $t=T/2$ are well represented with the proposed TS model when considering higher harmonic components. The relative difference between the solutions at $t=T/8$ is reduced to less than 1\% with $n=3$ (Fig.~\ref{InteriorFieldNolinearError}).   }

\begin{figure} [t]
	\hspace{-0.3cm}
	\includegraphics[trim={0cm 0cm 0cm 0cm},width=0.5\textwidth]{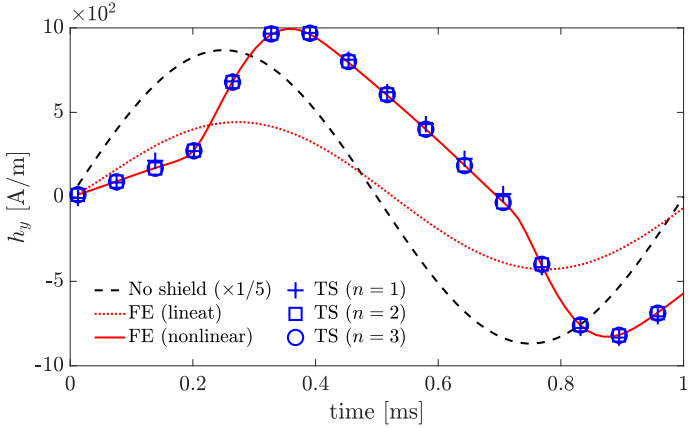}
	\caption{\small Time-evolution of $h_y$ at point $P_1$ for the nonlinear case. The solutions obtained  from the TS model with $n=1$ to 3 are compared with the 2-D FE solution. Solutions without the shield and in linear case ($\mu_r=$ 12500) are presented for the sake of comparison.}
	\label{ExteriorFieldNolinear}
	\vspace{-3mm}
\end{figure}

The time-evolution of $h_y$ at point $P_1$ is shown in Fig.~\ref{ExteriorFieldNolinear}. The solution of the 2-D FE problem without the shield gives the field at this point, which has the same waveform as the current~$I$ applied to the wires. The solution for a linear shield problem with $\mu_r=12500$ is presented for the sake of comparison with the nonlinear solution. Finally, the solution for the shield modeled with the TS model shows improvement as~$n$ increases, while high accuracy is observed when compared to the nonlinear reference solution. The saturation effect is clearly observed. 

\begin{figure} [t]
	\includegraphics[trim={0cm 0cm 0cm 0cm},width=0.49\textwidth]{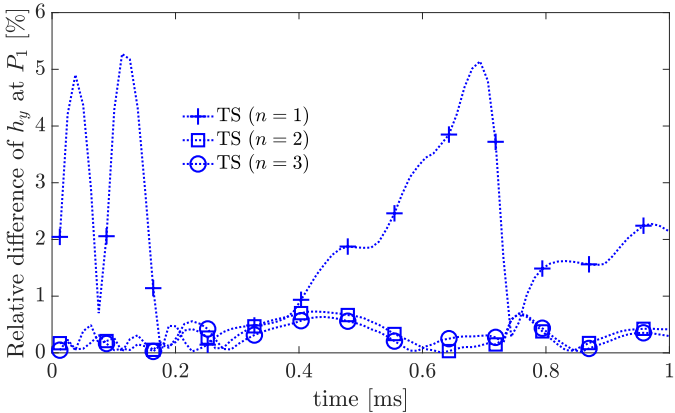}
	\caption{\small \textcolor{black}{Relative difference between the TS and the reference FE solutions for the time-evolution of $h_y$ at point $P_1$ for the nonlinear case.}}
	\label{ExteriorFieldNolinearError}
	\vspace{-3mm}
\end{figure}

The relative difference of the TS model to the FE solution at $P_1$ \textcolor{black}{at each time step} is presented in Fig.~\ref{ExteriorFieldNolinearError}. \textcolor{black}{It was reduced from more than $5\%$ with \mbox{$n=1$} to less than $1\%$ with $n\geqslant 2$.} With a suitable choice of the basis functions, the proposed TS model can certainly be an option for simulating nonlinear thin sheets.

\section{Conclusion}

In this paper, a time-domain extension of the classical TS model for thin sheets was elaborated and discussed using the $\bm{h}$-formulation. In our approach, the addition of~$n$ sets of basis functions derived from the steady-state solutions for the problem of a slab of finite thickness \marc{permits the representation of} the time evolution of the field quantities inside the thin region and its surroundings. We did apply this method to analyze the shielding efficiency of conducting and ferromagnetic planar sheets in harmonic and time-transient regimes for linear and nonlinear shield configurations.

In the harmonic regime, the proposed model gives~IBCs equivalents to those used in the classical TS model, since both methods are based on the solutions~\eqref{Sol1D_h} and~\eqref{Sol1D_e} and are directly included in the FE formulation. Our model, however, can also be used in time-transient FE analysis. We showed that, by adding a small number of hyperbolic basis functions, high precision can be achieved. \textcolor{black}{In all cases, the relative differences were $<3\%$ with the reference solutions with $n \geqslant3$, including the critical region near the extremities of the shield.} Furthermore, the proposed model \marc{can achieve comparable errors with less degrees of freedom, and hence at a lower computational cost, while also avoiding meshes with poor aspect ratios.} 

The TS model presented in this paper is still application-dependent since the set of hyperbolic basis functions must be defined according to the frequency content of the magnetic field inside the thin region, as well as the material composition of the sheet. 
\textcolor{black}{Although we did not yet find a general rule to select the basis functions, the latter are easy to derive since they originate from the analytic solutions of the 1-D linear flux diffusion problem in harmonic regime. Therefore, as long as the set of basis functions is rich enough to represent a diversity of penetration depths (which varies dynamically with local magnetic saturation), nonlinear solutions can be well approximated with this approach.}
Finally, the presented methodology can be easily extended to other FE formulations, such as the \mbox{$\bm{a}$-formulation,} as well as 3-D shielding problems.

\section*{Acknowledgment}

The authors would like to thank Prof. Christophe Geuzaine for fruitful discussions and for putting the Gmsh and GetDP codes in the public domain.

This work has been supported in part by the Coordenação de Aperfeiçoamento de Pessoal de Nível Superior – Brazil (CAPES) -  Finance code 001, and in part by the Fonds de Recherche du Québec - Nature et Technologies (FRQNT). Also, the collaboration between the authors was greatly facilitated by the MITACS Globalink internship program.

\ifCLASSOPTIONcaptionsoff
  \newpage
\fi

\bibliographystyle{IEEEtran}
\bibliography{IEEEabrv,Bibli}

\begin{thebibliography}{10}
\providecommand{\url}[1]{#1}
\csname url@samestyle\endcsname
\providecommand{\newblock}{\relax}
\providecommand{\bibinfo}[2]{#2}
\providecommand{\BIBentrySTDinterwordspacing}{\spaceskip=0pt\relax}
\providecommand{\BIBentryALTinterwordstretchfactor}{4}
\providecommand{\BIBentryALTinterwordspacing}{\spaceskip=\fontdimen2\font plus
\BIBentryALTinterwordstretchfactor\fontdimen3\font minus
  \fontdimen4\font\relax}
\providecommand{\BIBforeignlanguage}[2]{{%
\expandafter\ifx\csname l@#1\endcsname\relax
\typeout{** WARNING: IEEEtran.bst: No hyphenation pattern has been}%
\typeout{** loaded for the language `#1'. Using the pattern for}%
\typeout{** the default language instead.}%
\else
\language=\csname l@#1\endcsname
\fi
#2}}
\providecommand{\BIBdecl}{\relax}
\BIBdecl

\bibitem{bottauscio2006transient}
O.~Bottauscio, M.~Chiampi, and A.~Manzin, ``Transient analysis of thin layers
  for the magnetic field shielding,'' \emph{IEEE {T}ransactions on
  {M}agnetics}, vol.~42, no.~4, pp. 871--874, Mar 2006.

\bibitem{Igarashi1998a}
H.~Igarashi, A.~Kost, and T.~Honma, ``{A three dimensional analysis of magnetic
  fields around a thin magnetic conductive layer using vector potential},''
  \emph{IEEE {T}ransactions on {M}agnetics}, vol.~34, no.~5, pp. 2539--2542,
  1998.

\bibitem{Igarashi1998b}
------, ``{Impedance boundary condition for vector potentials on thin layers
  and its application to integral equations},'' \emph{EPJ Applied Physics},
  vol.~1, no.~1, pp. 103--109, 1998.

\bibitem{bottauscio2004numerical}
O.~Bottauscio, M.~Chiampi, and A.~Manzin, ``Numerical analysis of magnetic
  shielding efficiency of multilayered screens,'' \emph{IEEE {T}ransactions on
  {M}agnetics}, vol.~40, no.~2, pp. 726--729, 2004.

\bibitem{Rasilo2020}
P.~{Rasilo}, J.~{Vesa}, and J.~{Gyselinck}, ``Electromagnetic modeling of
  ferrites using shell elements and random grain structures,'' \emph{IEEE
  Transactions on Magnetics}, vol.~56, no.~2, pp. 1--4, 2020.

\bibitem{geuzaine2000dual}
C.~Geuzaine, P.~Dular, and W.~Legros, ``Dual formulations for the modeling of
  thin electromagnetic shells using edge elements,'' \emph{IEEE {T}ransactions
  on {M}agnetics}, vol.~36, no.~4, pp. 799--803, Jul 2000.

\bibitem{marchandise2014optimal}
E.~Marchandise, J.-F. Remacle, and C.~Geuzaine, ``Optimal parametrizations for
  surface remeshing,'' \emph{Engineering with Computers}, vol.~30, no.~3, pp.
  383--402, 2014.

\bibitem{krahenbuhl1993thin}
L.~Kr{\"a}henb{\"u}hl and D.~Muller, ``Thin layers in electrical engineering -
  example of shell models in analysing eddy-currents by boundary and finite
  element methods,'' \emph{IEEE {T}ransactions on {M}agnetics}, vol.~29, no.~2,
  pp. 1450--1455, Mar 1993.

\bibitem{mayergoyz1995calculation}
I.~D. Mayergoyz and G.~Bedrosian, ``On calculation of 3-{D} eddy currents in
  conducting and magnetic shells,'' \emph{IEEE {T}ransactions on {M}agnetics},
  vol.~31, no.~3, pp. 1319--1324, May 1995.

\bibitem{Guerin1994}
C.~Gu{\'e}rin, ``D{\'e}termination des pertes par courants de {F}oucault dans
  les cuves de transformateurs. mod{\'e}lisation de r{\'e}gions minces et prise
  en compte de la saturation des mat{\'e}riaux magn{\'e}tiques en r{\'e}gime
  harmonique,'' Ph.D. dissertation, Institut National Polytechnique de
  Grenoble-INPG, 1994.

\bibitem{Biro1997}
O.~B\'ir\'o, I.~Bardi, K.~Preis, W.~Renhart, and K.~Richter, ``A finite element
  formulation for eddy current carrying ferromagnetic thin sheets,'' \emph{IEEE
  Transactions on Magnetics}, vol.~33, no.~2, pp. 1173--1178, 1997.

\bibitem{gyselinck2008}
J.~Gyselinck, R.~V. Sabariego, P.~Dular, and C.~Geuzaine, ``Time-domain
  finite-element modeling of thin electromagnetic shells,'' \emph{IEEE
  {T}ransactions on {M}agnetics}, vol.~44, no.~6, pp. 742--745, June 2008.

\bibitem{sabariego2008h}
R.~V. Sabariego, C.~Geuzaine, P.~Dular, and J.~Gyselinck,
  ``\BIBforeignlanguage{English}{$h$- and $a$-formulations for the time-domain
  modelling of thin electromagnetic shells},''
  \emph{\BIBforeignlanguage{English}{IET Science, Measurement and Technology}},
  vol.~2, pp. 402--408, Nov 2008.

\bibitem{sabariego2009}
------, ``Nonlinear time-domain finite-element modeling of thin electromagnetic
  shells,'' \emph{IEEE {T}ransactions on {M}agnetics}, vol.~45, no.~3, pp.
  976--979, Mar 2009.

\bibitem{gyselinck2004time}
J.~Gyselinck and P.~Dular, ``A time-domain homogenization technique for
  laminated iron cores in 3-{D} finite-element models,'' \emph{IEEE
  {T}ransactions on {M}agnetics}, vol.~40, no.~2, pp. 856--859, Mar 2004.

\bibitem{gyselinck2006}
J.~Gyselinck, R.~V. Sabariego, and P.~Dular, ``A nonlinear time-domain
  homogenization technique for laminated iron cores in three-dimensional
  finite-element models,'' \emph{IEEE {T}ransactions on {M}agnetics}, vol.~42,
  no.~4, pp. 763--766, Apr 2006.

\bibitem{gyselinck2015finite}
J.~Gyselinck, P.~Dular, L.~Kr{\"a}henb{\"u}hl, and R.~V. Sabariego,
  ``Finite-element homogenization of laminated iron cores with inclusion of net
  circulating currents due to imperfect insulation,'' \emph{IEEE {T}ransactions
  on {M}agnetics}, vol.~52, no.~3, pp. 1--4, 2015.

\bibitem{brandt1996superconductors}
E.~H. Brandt, ``Superconductors of finite thickness in a perpendicular magnetic
  field: Strips and slabs,'' \emph{Physical review B}, vol.~54, no.~6, p. 4246,
  1996.

\bibitem{Zhang2017}
H.~Zhang, M.~Zhang, and W.~Yuan, ``An efficient 3d finite element method model
  based on the t--a formulation for superconducting coated conductors,''
  \emph{Superconductor Science and Technology}, vol.~30, no.~2, p. 024005,
  2017.

\bibitem{Liang2017}
F.~Liang, S.~Venuturumilli, H.~Zhang, M.~Zhang, J.~Kvitkovic, S.~Pamidi,
  Y.~Wang, and W.~Yuan, ``A finite element model for simulating second
  generation high temperature superconducting coils/stacks with large number of
  turns,'' \emph{Journal of Applied Physics}, vol. 122, no.~4, p. 043903, 2017.

\bibitem{Berrospe2019}
E.~Berrospe-Juarez, V.~M. Zerme{\~n}o, F.~Trillaud, and F.~Grilli, ``Real-time
  simulation of large-scale hts systems: multi-scale and homogeneous models
  using the t--a formulation,'' \emph{Superconductor Science and Technology},
  vol.~32, no.~6, p. 065003, 2019.

\bibitem{Sabariego2010}
R.~V. {Sabariego}, P.~{Dular}, C.~{Geuzaine}, and J.~{Gyselinck},
  ``Surface-impedance boundary conditions in dual time-domain finite-element
  formulations,'' \emph{IEEE {T}ransactions on {M}agnetics}, vol.~46, no.~8,
  pp. 3524--3531, 2010.

\bibitem{mayergoyz1994finite}
I.~Mayergoyz and G.~Bedrosian, ``On finite element implementation of impedance
  boundary conditions,'' \emph{Journal of Applied Physics}, vol.~75, no.~10,
  pp. 6027--6029, 1994.

\bibitem{knoepfel2008magnetic}
H.~E. Knoepfel, \emph{Magnetic fields: a comprehensive theoretical treatise for
  practical use}.\hskip 1em plus 0.5em minus 0.4em\relax John Wiley \& Sons,
  2008.

\bibitem{Gyselinck2009}
J.~{Gyselinck}, P.~{Dular}, C.~{Geuzaine}, and R.~V. {Sabariego},
  ``Surface-impedance boundary conditions in time-domain finite-element
  calculations using the magnetic-vector-potential formulation,'' \emph{IEEE
  Transactions on Magnetics}, vol.~45, no.~3, pp. 1280--1283, 2009.

\bibitem{Dular2011}
P.~Dular, V.~Q. Dang, R.~V. Sabariego, L.~Kr{\"a}henb{\"u}hl, and C.~Geuzaine,
  ``Correction of thin shell finite element magnetic models via a subproblem
  method,'' \emph{IEEE {T}ransactions on {M}agnetics}, vol.~47, no.~5, pp.
  1158--1161, 2011.

\bibitem{Gmsh2009}
C.~Geuzaine and J.-F. Remacle, ``Gmsh: A 3-{D} finite element mesh generator
  with built-in pre-and post-processing facilities,'' \emph{International
  journal for numerical methods in engineering}, vol.~79, no.~11, pp.
  1309--1331, 2009.

\bibitem{GetDP2013}
P.~Dular and C.~Geuzaine, ``Get{D}{P} reference manual: the documentation for
  {G}et{D}{P}, a general environment for the treatment of discrete problems,''
  \emph{University of Li\`ege}, 2013.

\bibitem{Dular2020}
J.~{Dular}, C.~{Geuzaine}, and B.~{Vanderheyden}, ``Finite-element formulations
  for systems with high-temperature superconductors,'' \emph{IEEE Transactions
  on Applied Superconductivity}, vol.~30, no.~3, pp. 1--13, 2020.

\end{thebibliography}

\end{document}